\documentclass[12pt]{article}
\usepackage{amssymb,amsmath,amsopn,amsfonts,mathrsfs,bm}
\usepackage[dvips]{color}
\usepackage{colordvi,multicol}
\usepackage{epsf}
\usepackage{hyperref}
\usepackage{bbm}
 \usepackage[titletoc]{appendix}
 \usepackage{color}

\newfam\msbfam
\font\tenmsb=msbm10 \textfont\msbfam=\tenmsb \font\sevenmsb=msbm7
\scriptfont\msbfam=\sevenmsb \font\fivemsb=msbm5
\scriptscriptfont\msbfam=\fivemsb

\newcommand{\normmm}[1]{{\left\vert\kern-0.25ex\left\vert\kern-0.25ex\left\vert #1 
    \right\vert\kern-0.25ex\right\vert\kern-0.25ex\right\vert}}

\def\proof{\vspace{1mm}\noindent{\it Proof}\quad}

 \newtheorem{lemma}{Lemma}[section]
 \newtheorem{proposition}[lemma]{Proposition}
 \newtheorem{theorem}[lemma]{Theorem}
 \newtheorem{corollary}[lemma]{Corollary}
 \newtheorem{remark}[lemma]{Remark}
 \newtheorem{definition}[lemma]{Definition}

\numberwithin{equation}{section}

\def\bc{\begin{center}}
\def\ec{\end{center}}
\def\no{\noindent}

  \makeatletter
\newcommand{\rmnum}[1]{\romannumeral #1}
\newcommand{\Rmnum}[1]{\expandafter\@slowromancap\romannumeral #1@}
\makeatother

\begin{document}
\bigbreak
\title {{\bf Weak solutions to the sharp interface limit of stochastic Cahn-Hilliard equations
\thanks{Research supported in part  by NSFC (No.11671035). Financial support by the DFG through the CRC 1283 "Taming uncertainty and profiting from randomness and low regularity in analysis, stochastics and their applications" is acknowledged.}\\}}
\author{{\bf Huanyu Yang}$^{\mbox{a,c,d}}$,{\bf Rongchan Zhu}$^{\mbox{b,c}}$,
\date {}
\thanks{E-mail address: hyang@math.uni-bielefeld.de(H. Y. Yang), zhurongchan@126.com(R. C. Zhu),}\\ \\
\small $^{\mbox{a}}$School of Mathematical Science, University of Chinese Academy of Sciences, Beijing 100049, China,\\
\small $^{\mbox{b}}$Department of Mathematics, Beijing Institute of Technology, Beijing 100081, China,\\
\small $^{\mbox{c}}$ Department of Mathematics, University of Bielefeld, D-33615 Bielefeld, Germany,\\
\small $^{\mbox{d}}$ Academy of Mathematics and Systems Science, Chinese Academy of Sciences, Beijing 100190, China,}

\maketitle
\noindent {\bf Abstract}

We study the asymptotic limit, as $\varepsilon\searrow 0$, of solutions of the stochastic Cahn-Hilliard equation:
$$
   \partial_t u^\varepsilon=\Delta \left(-\varepsilon\Delta u^\varepsilon+\frac{1}{\varepsilon}f(u^\varepsilon)\right)+\dot{\mathcal{W}}^\varepsilon_t, \\
$$
where $\mathcal{W}^\varepsilon=\varepsilon^\sigma W$ or $\mathcal{W}^\varepsilon=\varepsilon^\sigma W^\varepsilon$, $W$ is a $Q$-Wiener process and $W^\varepsilon$ is smooth in time and converges to $W$ as $\varepsilon\searrow 0$. In the case that $\mathcal{W}^\varepsilon=\varepsilon^\sigma W$, we prove that for all $\sigma>\frac{1}{2}$, the solution $u^\varepsilon$ converges to a weak solution to an appropriately defined limit of the deterministic Cahn-Hilliard equation. In radial symmetric case we prove that for all $\sigma\geq\frac{1}{2}$, $u^\varepsilon$ converges to the deterministic Hele-Shaw model. In the case that $\mathcal{W}^\varepsilon=\varepsilon^\sigma W^\varepsilon$, we prove that for all $\sigma>0$, $u^\varepsilon$ converges to the weak solution to the deterministic limit Cahn-Hilliard equation. In radial symmetric case we prove that $u^\varepsilon$ converges to deterministic Hele-Shaw model when $\sigma>0$ and converges to  a stochastic model related to stochastic Hele-Shaw model when $\sigma=0$.

 \no{\footnotesize{\bf Keywords}:\hspace{2mm} the deterministic/stochastic Cahn-Hilliard equation, sharp interface limit, Hele-Shaw model, tightness, varifold

\normalsize
\section{Introduction}\label{s1}
We consider the sharp interface limit of the following stochastic Cahn-Hilliard equation on a bounded smooth open domain $\mathcal{D}\subset\mathbb{R}^d$ ($d=2,3$):
\begin{equation}\label{1.1}
\left\{
\begin{aligned}
&du^{\varepsilon}=\Delta v^\varepsilon dt+\varepsilon^\sigma dW_t, \quad (t,x)\in[0,T]\times\mathcal{D},\\
&v^\varepsilon=-\varepsilon\Delta u^{\varepsilon}(t)+\frac{1}{\varepsilon} f(u^{\varepsilon}(t)),\quad (t,x)\in[0,T]\times\mathcal{D},\\
&\frac{\partial u^{\varepsilon}}{\partial n}=\frac{\partial v^{\varepsilon}}{\partial n}=0, \quad (t,x)\in[0,T]\times\partial\mathcal{D},\\
&u^{\varepsilon}(0,x)=u_0^\varepsilon(x),\quad x\in\mathcal{D}.
\end{aligned}\right.
\end{equation}
Here $W$ is a $Q$-Wiener process where $Q$ satisfies (\ref{2.1}) and (\ref{2.2}). 
$f(u)=F'(u)$ where $F(u)=\frac{1}{4}(u^2-1)^2$ is the double well potential and the initial data $u_0^\varepsilon$ satisfies
\begin{equation}\label{1.2}
\left\{
\begin{aligned}
&\sup_{0<\varepsilon\leq 1}\int_\mathcal{D}\left(\frac{\varepsilon}{2}|\nabla u_0^\varepsilon(x)|^2+\frac{1}{\varepsilon}F(u_0^\varepsilon(x))\right)dx\leq\mathcal{E}_0<\infty,\\
&\frac{1}{|\mathcal{D}|}\int_\mathcal{D}u_0^\varepsilon(x)dx=m_0\in(-1,1)\quad\forall\varepsilon\in(0,1].
\end{aligned}
\right.
\end{equation}

\vskip.10in
\textbf{Sharp interface limit of deterministic equation.}
The deterministic Cahn-Hilliard equation
\begin{equation}\label{1.4}
\left\{
\begin{aligned}
&\partial_t u^{\varepsilon}=\Delta v^\varepsilon, \quad (t,x)\in[0,T]\times\mathcal{D},\\
&v^\varepsilon=-\varepsilon\Delta u^{\varepsilon}(t)+\frac{1}{\varepsilon} f(u^{\varepsilon}(t)),\quad (t,x)\in[0,T]\times\mathcal{D},\\
&\frac{\partial u^{\varepsilon}}{\partial n}=\frac{\partial v^{\varepsilon}}{\partial n}=0, \quad (t,x)\in[0,T]\times\partial\mathcal{D},\\
&u^{\varepsilon}(0,x)=u_0^\varepsilon(x),\quad x\in\mathcal{D},
\end{aligned}\right.
\end{equation}
 is widely accepted as a good model to describe the complicated phase separation and coarsening phenomena in a melted alloy that is quenched to a temperature at which only two different concentration phases can exist stably. It was formally derived by Pego \cite{London:wf} and rigorous proved by \cite{ALIKAKOS:1994vc} by using the method of matched asymptotic expansions that the equation (\ref{1.4}) converges to the motion by Mullins-Sekerka law. That is, as $\varepsilon\searrow 0$, the chemical potential $v^\varepsilon$ tends to a limit $v$ which, together with a free boundary $\Gamma:=\cup_{0\leq t\leq T}(\{t\}\times\Gamma_t)$, solves the following free boundary problem:
\begin{equation}\label{1.3}
\left\{
   \begin{aligned}
   \Delta v&=0 \;\text{in}\;\mathcal{D}\setminus\Gamma_t,\; t>0,\\
\frac{\partial v}{\partial n}&=0\;\text{on}\;\partial\mathcal{D},\\
   v&=SH\;\text{on}\;\Gamma_t,\\
  \mathcal{V}&=\frac{1}{2}(\partial_{n}v^+-\partial_{n}v^-)\;\text{on}\;\Gamma_t,\\
   \end{aligned}
   \right.
\end{equation}
where 
$$
S=\int_{-1}^1\sqrt{\frac{F(s)}{2}}ds=\frac{2}{3},
$$
$H$ is the scalar mean curvature of $\Gamma_t$ with the sign convention that convex hypersurfaces have positive mean curvature, $\mathcal{V}$ is the normal velocity of the interface with the sign convention that the normal velocity of expanding hypersurfaces is positive, $n$ is the unit ourward normal either to $\partial\mathcal{D}$ or to $\Gamma_t$. Denote $\mathcal{D}^+$ be the region enclosed by $\Gamma_t$ and $\mathcal{D}^-=\mathcal{D}\setminus(\mathcal{D}^+\cup\Gamma_t)$ is the region  and  $v^+$, $v^-$ are respectively the restriction of $v$ on $[0,t]\times\mathcal{D}^+$ and $[0,t]\times\mathcal{D}^-$.

Later in \cite{Chen:1996bg}, the author formulate a weak solution to the free boundary problem (\ref{1.3}) (see Definition \ref{d2.1}) and show that the solutions of (\ref{1.4}) approach, as $\varepsilon\searrow 0$, to weak solutions of (\ref{1.3}) by using a compactness argument. In fact, the Cahn-Hilliard equation (\ref{1.4}) is an $H^{-1}$-gradient flow with the van der Waals-Cahn-Hilliard energy functional
\begin{equation}\label{1.5}
\mathcal{E}^\varepsilon(u^\varepsilon):=\int_\mathcal{D}e^\varepsilon(u^\varepsilon)dx,\quad e^\varepsilon(u^\varepsilon):=\frac{\varepsilon}{2}|\nabla u^\varepsilon|^2+\frac{1}{\varepsilon}F(u^\varepsilon).
\end{equation}
Denote by $(u_D^\varepsilon,v_D^\varepsilon)$ the solution to the deterministic Cahn-Hilliard equation (\ref{1.4}).
One can directly verify that 
\begin{equation}
\frac{d}{dt}\mathcal{E}^\varepsilon(u_D^\varepsilon)=-\int_\mathcal{D}|\nabla v_D^\varepsilon|^2\leq 0,
\end{equation} 
which is also called the Lyapunov property for equation (\ref{1.4}). Thus $\mathcal{E}^\varepsilon(u_D^\varepsilon)$ is uniformly bounded in $t,\varepsilon>0$ with the initial value satisfying (\ref{1.2}).  Note that as $\varepsilon\to 0$, $F(u_D^\varepsilon)\to 0$, which is equivalent to $u_D^\varepsilon\to-1+2\chi_{E}$ for some $E\subset[0,T]\times\mathcal{D}$ where $\chi_{E}$ is the characteristic function of $E$, i.e. $\chi_{E}(x)=1$ when $x\in E$ and $\chi_{E}(x)=0$ when $x\not\in E$. $\Gamma_t:=\partial E_t$ is the interface. By using a varifold approach, Chen in \cite{Chen:1996bg} analyzed the property of the limit of the solutions to equation (\ref{1.4}) and then proposed a definition of weak solution of this limit. Any classical smooth solutions to (\ref{1.4}) are weak solutions. In some special case, the smooth weak solutions are also classical solutions to (\ref{1.3}).
We need to mention that in \cite{ALIKAKOS:1994vc}, the convergence of solutions to Cahn-Hilliard equation (\ref{1.4}) to (\ref{1.3}) is proved under the assumption on the existence of smooth solution to  (\ref{1.4}). While in \cite{Chen:1996bg}, Chen proved the convergence of the solution to equation (\ref{1.4}) and analyzed the limit directly. No assumption on existence of solution to  (\ref{1.4}) is required in \cite{Chen:1996bg}. 
\vskip.10in
\textbf{Sharp interface limit of stochastic equation.}
For stochastic case,   the sharp interface limit of equation (\ref{1.1}) was first considered in \cite{Antonopoulou:2018gh}, where the authors compare the solutions to equation (\ref{1.1}) with the approximation solutions constructed in \cite{ALIKAKOS:1994vc}. They proved that if the smooth solution to (\ref{1.3}) exists, then for all $\sigma>\frac{23}{3}$ when $d=2$ and all $\sigma>11$ when $d=3$, the solutions to (\ref{1.1}) will converge to the solutions to (\ref{1.3}). Later in \cite{Banas:2019wz}, the authors extend their results to the case that $W$ is a cylindrical Wiener process or conservative noise for $\sigma$ large enough. For $\sigma\geq0$ small, the perturbation by the noise become much stronger. It is reasonable to think that the equation to (\ref{1.1}) will no longer converge to (\ref{1.3}) when $\sigma$ is small.
But the method in\cite{Antonopoulou:2018gh} can be only applied to prove the convergence to  (\ref{1.3}) and also seems not easy to obtain the convergence for small $\sigma$.

In this paper, for small $\sigma\geq0$, we extend the method in \cite{Chen:1996bg} to equation (\ref{1.1}) and obtain weak solutions to the limit of equation (\ref{1.1}). Then we consider the limit of the solution to equation (\ref{1.1}) directly, which enables us to analyze different models  the limit should satisfy.
We mainly consider (\ref{1.1} with two types driven noise: $Q$-Wiener process and "smeared" noise which is smooth in time.
\vskip.10in
\textbf{The equation with $Q$-Wiener process for \bm{$\sigma\geq\frac{1}{2}$}.} In this case, we can obtain that for $\sigma>\frac{1}{2}$, the solutions to equation (\ref{1.1}) converge to the weak solutions defined in Definition \ref{d2.1}.  In fact,
motivated by \cite{DaPrato:1996kk}, we apply the It\^o's formula to $\mathcal{E}^\varepsilon(u^\varepsilon)$ and prove the Lyapunov property of equation (\ref{1.1}) for all $\bm{\sigma\geq\frac{1}{2}}$ (see Lemma \ref{l3.1}). By a tightness argument, we prove that for all $\sigma>\frac{1}{2}$, the solutions to equation (\ref{1.1}) converge to the weak solution of the limit of deterministic Cahn-Hilliard equation (\ref{1.4}) defined by Chen \cite{Chen:1996bg} (see Theorem \ref{t2.3}). For $\sigma=\frac{1}{2}$, the tightness and convergence results are still true. But we cannot conclude that the limit is a weak solution defined in Definition \ref{d2.1}.

Particularly in radial symmetric case, we prove that for all $\sigma\geq\frac{1}{2}$, the limit of solutions to equation (\ref{1.1}) satisfy (\ref{1.4}) in the weak sense. 
Thus we conjecture that in general
for $\mathbb{P}-a.s.\;\omega$,  as $\varepsilon\searrow 0$, the chemical potential $v^\varepsilon(\omega)$ tends to a limit $v(\omega)$ which, together with a free boundary $\Gamma(\omega):=\cup_{0\leq t\leq T}(\{t\}\times\Gamma_t(\omega))$,
$(v(\omega),\Gamma(\omega))$ satisfies (\ref{1.4}). 

\vskip.10in

\textbf{The equation with "smeared" noise for \bm{$\sigma\geq 0$}.}
Moreover, we consider stochastic Cahn-Hilliard equation driven by "smeared" noise which is smooth in time. This kind of noise were considered also for stochastic Allen-Cahn equation in \cite{FUNAKI:1999ho,Weber:2010dl,Funaki:2019is}. 

We smoothen the noise in time and consider the following random PDE:
\begin{equation}\label{1.8}
\left\{
\begin{aligned}
&\frac{\partial u^{\varepsilon}}{\partial t}=\Delta v^\varepsilon+\varepsilon^\sigma \xi^\varepsilon_t, \quad (t,x)\in[0,T]\times\mathcal{D},\\
&v^\varepsilon=-\varepsilon\Delta u^{\varepsilon}(t)+\frac{1}{\varepsilon} f(u^{\varepsilon}(t)),\quad (t,x)\in[0,T]\times\mathcal{D},\\
&\frac{\partial u^{\varepsilon}}{\partial n}=\frac{\partial v^{\varepsilon}}{\partial n}=0, \quad (t,x)\in[0,T]\times\partial\mathcal{D},\\
&u^{\varepsilon}(0,x)=u_0^\varepsilon(x),\quad x\in\mathcal{D},
\end{aligned}\right.
\end{equation}
where $\xi^\varepsilon_t=\frac{dW^\varepsilon}{dt}$,  $W^\varepsilon_t:=\int_{-\infty}^\infty\rho_\varepsilon(t-s)W_sds$ and $\rho_\varepsilon$ is an approximate delta function on $\mathbb{R}$. Formally as $\varepsilon\to 0$, $\xi^\varepsilon\to \frac{dW}{dt}$.
Since $\xi^\varepsilon_t$ is smooth in time, this enables us to apply the Newton-Leibniz formula to $\mathcal{E}^\varepsilon(u^\varepsilon)$ and obtain the Lyapunov property. Thus the tightness and the convergence results hold for all $\bm{\sigma\geq 0}$. Similar as before, for all $\sigma>0$, the solutions to (\ref{1.8})  converge to the weak solution to Definition \ref{d2.1} (see Theorem \ref{t6.3}).
For the interesting case that \bm{$\sigma=0$}, when $\varepsilon\searrow 0$, we have that $u^\varepsilon\to -1+2\chi_E$ for some $E\in[0,T]\times\mathcal{D}$, $v^\varepsilon\to v$ and 
\begin{equation}\label{1.7}
2d\chi_E=\Delta vdt+dW_t.
\end{equation}
(\ref{1.7}) actually gives a weak formula to describe how the evolution of the interface $\Gamma_t:=\partial E_t$ is governed by the noise $W$ (see Theorem \ref{t6.4}). This gives the first rigorous result of the sharp interface limit of stochastic Cahn-Hillliard limit to a stochastic model. Similar as before, we  conjecture that for $\mathbb{P}-a.s.\;\omega$,  as $\varepsilon\searrow 0$, the chemical potential $v^\varepsilon(\omega)$ tends to a limit $v(\omega)$ which, together with a free boundary $\Gamma(\omega):=\cup_{0\leq t\leq T}(\{t\}\times\Gamma_t(\omega))$,
$(v(\omega),\Gamma(\omega))$ satisfies the following stochastic problem:
\begin{equation}\label{1.9}
\left\{
   \begin{aligned}
   \Delta vdt&=-dW_t \;\text{in}\;\mathcal{D}\setminus\Gamma_t,\; t>0,\\
\frac{\partial v}{\partial n}&=0\;\text{on}\;\partial\mathcal{D},\\
   v&=SH\;\text{on}\;\Gamma_t,\\
 \mathcal{V}dt&=\frac{1}{2}\left[\frac{\partial}{\partial n}\right]_{\Gamma_t}(vdt+\Delta^{-1} d{W}_t).
   \end{aligned}
   \right.
\end{equation}
\vskip.10in
\vskip.10in
We also mention that Chen's definition for weak solution in Definition \ref{d2.1} is not so ``perfect", since it is still unknown whether in general such a smooth weak solution is a classical solution to (\ref{1.3}). The problems come from that a "good" weak formula for the third equation in (\ref{1.3}) is still missing. 
Moreover, in \cite{Antonopoulou:2018gh} the authors also give some different conjectures about the sharp interface limit of equation (\ref{1.1}) via a formal calculation, especially in the case that $\sigma=1$. In their case the value of $v$ on the interface is different from ours. As what we analyze in Remark \ref{r5.9}, our model (\ref{1.9}) fit quite well in radial symmetric case. But in general case, we still cannot give a fully rigorous proof.

In fact, identifying the value of $v$ on the interface $\Gamma_t$ is the main task of varifold approach to study the sharp interface limit of both Cahn-Hilliard equation and Allen-Cahn equation (cf. \cite{Hutchinson:2000vy,Tonegawa:2002dl,Tonegawa:2005ut,Roger:2006hz,Le:2008dt,Roger:2008du}). In these literature, the authors give a weak formula for the third equation in (\ref{1.3}). But they are unable to prove the limit of the solutions to  equation (\ref{1.4}) satisfy such weak formula. Until now, a fully rigorous proof of the sharp interface limit of Cahn-Hilliard equation is still missing.

The problem (\ref{1.3}) is often called the Mullins-Sekera problem or the two-phase Hele-Shaw problem. The local existence and uniqueness of classical solutions of (\ref{1.3}) when the initial hypersufrace $\Gamma_0$ is sufficiently smooth can be found in \cite{Xinfu:1996ci,Escher:1997vw}. For general initial hypersurfaces $\Gamma_0$, existence of weak solutions of (\ref{1.3}) was proved in \cite{Chen:1996bg,Roger:2005dw} by using the tools of varifolds. For the stochastic Hele-Shaw model (\ref{1.9}), there is no result for existence.

Finally, as what we mentioned before, the methods in \cite{Chen:1996bg} and also in this paper are deeply related to the theory of varifolds. We recall some related definitions in Section \ref{s2}. In fact,
varifolds represent very natural generalizations of classical $n$-surfaces, as they encode, loosely speaking, a joint distribution of mass and tangents. More technically, varifolds are Radon measures defined on the Grassmann bundle $\mathbb{R}^d\times G({n,d})$, whose elements are pairs $(x, S)$ specifying a position in space and an unoriented $n$-plane. Varifolds have been proposed more than $50$ years ago by Almgren \cite{ALMGREN:1965td} as a mathematical model for soap films, bubble clusters, crystals, and grain boundaries. After Allard's fundamental work \cite{mathematics:vh}, varifolds have been successfully used in the context of Geometric Measure Theory, Geometric Analysis, and Calculus of Variations. One successful application of varifolds resulted in the definition and the study of a general weak mean curvature flow in \cite{Brakke:1978to}, which allowed to prove existence of mean curvature evolution with singularities in \cite{Kim:2017gz}. Beyond the theory of rectifiable varifolds, the flexibility of the varifold structure has been proved to be relevant to model diffuse interfaces, e.g., phase field approximations, and a crucial part in the proof of the convergence of the Allen-Cahn equation to Brakke’s mean curvature flow \cite{Ilmanen:1993wj,Tonegawa:2003hb,Takasao:2015fy}, or in the proof of the $\Gamma$-convergence of Cahn-Hilliard type energies to the Willmore energy (up to an additional perimeter term) \cite{Tonegawa:2005ut,Roger:2006hz,Le:2008dt,Roger:2008du}.

This paper is organized as follows: In Section \ref{s2} we give some basic notations and recall some definitions from geometric measure theory. In subsection \ref{ss2.3} we give the main results for (\ref{1.1}) driven by $Q$-Wiener process. In Section \ref{s3}, we establish certain $\varepsilon$-independent estimates for the solution to (\ref{1.1}), which allow us obtain tightness and then apply Skorohod's theorem to obtain a convergence subsequence for all $\sigma\geq\frac{1}{2}$. Moreover  for $\sigma>\frac{1}{2}$, we prove that this limit is actually a weak  solution to (\ref{1.3}). Similar as in \cite{Chen:1996bg}, in Section \ref{s4}, we study the radially symmetric case and prove that for all $\sigma\geq\frac{1}{2}$, the limit of the solution to equation (\ref{1.1}) satisfies the deterministic Hele-Shaw model (\ref{1.3}).
Finally in Section \ref{s6}, we consider the case for "smeared" noise $\xi^\varepsilon_t$ and obtain the convergence result for all $\sigma\geq 0$. For $\sigma>0$, the limit of the solution to (\ref{1.8}) is a weak solution to equation (\ref{1.3}). For $\sigma=0$, we obtain a stochastic characterisation of the evolution of the interface  (\ref{1.7}) and prove in radial symmetric case that it satisfies the stochastic Hele-Shaw model (\ref{1.9}).

\section{Preliminary}\label{s2}
\subsection{Basic notations and assumptions}

In the following, we denote by $S^{d-1}$ the unit sphere in $\mathbb{R}^d$ and $\vec{\nu}$ a generic element in $S^{d-1}$.  If $n=(n^1,\cdots,n^d)$, we denote by $n\otimes n$ the matrix $(n^in^j)_{d\times d}$. We use ``$\mathrm{I}$" to denote the identity matrix $(\delta_{ij})_{d\times d}$. For any $d\times d$ matrices $A=(a_{ij})$ and $B=(b_{ij})$, 
$$
A:B:=\text{Trace}(A^{T}B)=\sum_{i,j=1}^da_{ij}b_{ij}.
$$

We denote by $C_c^m(\mathcal{O})$ the space of $m$-th differentiable functions with compact support in $\mathcal{O}$ where $\mathcal{O}$ can be open or closed. 
Note that if $\mathcal{O}$ is compact, $C_c^m(\mathcal{O})=C^m(\mathcal{O})$. Moreover, we say a vector function $\vec{Y}=(Y^1,\cdots,Y^d)\in C_c^m(\mathcal{O};\mathbb{R}^d)$ if $Y^i\in C_c^m(\mathcal{O})$ for any $i=1,\cdots,d$.
For any $t>0$, we denote $\mathcal{O}_t:=[0,t]\times\mathcal{O}$.
We also denote by $\chi_E$ the characteristic function of a set $E$, which is defined by $\chi_E(x)=1$ for $x\in E$ and $
\chi_E(x)=0$ for $x\not\in E$.

Moreover, we denote $\mathcal{H}^n$ as the $n$-dimensional Hausdorff measure on $\mathbb{R}^d$ for any $n\in[0,d]$. For $d=n$, $\mathcal{H}^d$ is just the Lebesgue measure on $\mathbb{R}^d$.

We assume that $\mathcal{D}$ is a smooth bounded open domain in $\mathbb{R}^d$ $(d=2,3)$. Let $Q$ be an linear operator on $L^2(\mathcal{D})$, which is commuted with $\Delta$ and satisfies
\begin{equation}\label{2.1}
Q\mathbbm{1}=0,
\end{equation}
where $\mathbbm{1}(x)\equiv 1$ for any $x\in\mathcal{D}$
and
\begin{equation}\label{2.2}
\text{Tr}((-\Delta)Q)<+\infty.
\end{equation}
Let $(\Omega,\mathcal{F},\mathbb{P})$ be a stochastic basis and defined on it a $Q$-Wiener process $W$ on $L^2(\mathcal{D})$. 

According to \cite[Remark 2.2]{DaPrato:1996kk}, we have that 
\begin{theorem}Assume that $Q$ satisfies (\ref{2.1}), (\ref{2.2}), then for $\mathbb{P}-a.s.\;\omega$, the
equation (\ref{1.1}) has an analytic weak solution $u^{\varepsilon}\in C([0,T];H^1\cap L^4)$.
\end{theorem}

Let $u^\varepsilon$ be the solution to equation (\ref{1.1}), we set
\begin{equation}\label{3.1}
\mathcal{E}^\varepsilon(t):=\mathcal{E}^\varepsilon(u^\varepsilon)(t)=\int_\mathcal{D}e^\varepsilon(u^\varepsilon(t,x))dx,\quad e^\varepsilon(u^\varepsilon):=\frac{\varepsilon}{2}|\nabla u^\varepsilon|^2+\frac{1}{\varepsilon}F(u^\varepsilon).
\end{equation}

In the following, we recall several definitions from geometric measure theory (cf. \cite{Federer:2014iq,Simon:1983ti}).

 \vskip.10in
\textbf{Radon measures}

Let $\mathcal{O}$ be either an open or a closed domain. If $L$ is a bounded linear functional on $C_c(\mathcal{O})$ satisfying $\langle L,\psi\rangle\geq 0$ whenever $\psi\geq 0$ and $\psi\in C_c(\mathcal{O})$, the measure $\mu$ generated by 
$$
\mu(A)=\sup_{\psi\in C_c(A),|\psi|\leq 1}\langle L,\psi\rangle\quad \text{for all}\;\;A\;\;\text{open in}\;\;\mathcal{O}
$$
is called Radon measure on $\mathcal{O}$. We use $\langle\mu,\psi\rangle$ $\psi\in C_c(\mathcal{O})$ to denote the value $\int_{\mathcal{O}}\psi d\mu(=\langle L,\psi\rangle)$.

Let $\mathfrak{M}(\mathcal{D}_T)$ be the space of all finite signed measures on $\mathcal{D}_T$ and $\mathfrak{M}_R(\mathcal{D}_T)\subset\mathfrak{M}(\mathcal{D}_T)$ is the space of all Radon measures on $\mathcal{D}_T$. $\mathfrak{M}_R(\mathcal{D}_T)$ and $\mathfrak{M}(\mathcal{D}_T)$ are equipped with the total variation norm $\|\cdot\|_{TV}$ and weak topology, respectively.
\vskip.10in
\textbf{BV functions}

Let $u\in L^1(\mathcal{D})$. If the distributional gradient $Du$ defined by
$$
\langle Du,\vec{Y}\rangle:=\langle u,-\text{div}\vec{Y}\rangle\quad\forall\vec{Y}\in C_c^1(\mathcal{D};\mathbb{R}^d)
$$
can be extended as a bounded linear functional over $C_c(\mathcal{D};\mathbb{R}^d)$, then we say that $u$ is a function of bounded variation, denoted by $u\in BV(\mathcal{D})$. If $u\in BV(\mathcal{D})$, we use $D_iu$ to denote the measure on $C_c(\mathcal{D})$ generated by the functional $\langle u,-\partial_{x_i}\psi\rangle$ for all $\psi\in C_c^1(\mathcal{D})$. We denote by $|Du|$ the Radon measure generated by 
$$
|D(u)|(A):=\sup_{\vec{Y}\in C_c(A;\mathbb{R}^d),|\vec{Y}|\leq 1}\int_A u\;\text{div}\vec{Y}dx,\quad \forall A\;\;\text{open}\;\subset\mathcal{D}.
$$
One can show in \cite{Federer:2014iq} that $D_iu$ is absolutely continuous with respect to $|Du|$ and there exists a $|Du|$-measurable unit vector valued function $\vec{\nu}$ such that $Du=\vec{\nu}|Du|$, $|Du|-a.e.$.

We say that a set $E\subset\mathcal{D}$ is a BV set if $\chi_E\in BV(\mathcal{D})$. We denote $\vec{\nu}_E$ by 
\begin{equation}\label{normal}
D\chi_E=\vec{\nu}_E|D\chi_E|\;\;\text{or}\;\;\vec{\nu}_E(x):=\frac{D\chi_E(x)dx}{|D\chi_E|(x)dx}.
\end{equation}
Clearly, in the case that $\partial E$ is smooth, $\vec{\nu}_E$ is the unit inward normal of $E$ on $\partial E$.

In the following,
we introduce the several concepts of general varifold, which can be found in \cite[Chapter 8]{Simon:1983ti}.
\vskip.10in
\textbf{Varifolds}

Let $G(d,d-1)$ be the Grassmannian space which parametrizes of all $(d-1)$-dimensional linear subspaces of $\mathbb{R}^d$, which is a compact smooth manifold. For any $T\in G(d,d-1)$, $T$ can be identified with its unit normal vector $\vec{\nu}$. More precisely, $G(d,d-1)\cong P$, where $P:=S^{d-1}/\{\vec{\nu},-\vec{\nu}\}$ is the set of unit normals of unoriented $(d-1)$-planes in $\mathbb{R}^d$.

\begin{definition}\label{va}(varifold).
A varifold (or, more precisely a ($d-1$)-varifold) $V$ is a non-negative Radon measure on $G_{d-1}(\mathcal{D}):=\mathcal{D}\times G(d,d-1)$. The convergence of a sequence of varifolds is defined as the weak convergence in the sense of Radon measure.
\end{definition}

\begin{definition}\label{mass}(mass). Given a $(d-1)$-varifold $V$, there corresponds a Radon measure $\|V\|$ on $\mathcal{D}$ defined by 
 $$
\|V\|(A):=V(\pi^{-1}(A)), 
 $$
where $\pi$ is the projection $G_{d-1}(\mathcal{D})\ni(x,T)\longmapsto x$ onto $\mathcal{D}$. 
\end{definition}

\vskip.10in
\textbf{First variation of a varifold}
\begin{definition}\label{fv}
The first variation of a $(d-1)$-varifold $V$ in $\mathcal{D}$ is the vector-valued distribution $\delta V$ defined for any vector field $\vec{Y}=(Y^1,\cdots,Y^d)\in C_c^1(\mathcal{D},\mathbb{R}^d)$ as
$$
\langle\delta V,\vec{Y}\rangle:=\int_{G_{d-1}(\mathcal{D})}\mathrm{div}_T\vec{Y}(x)dV(x,T).
$$
Here for any $T\in G(d,d-1)$,
$$
\mathrm{div}_T\vec{Y}=\sum_{i=1}^d\nabla^T_iY^i,
$$
where $\nabla^T_i:=e_i\cdot\nabla^T$, $\{e_i\}_{i=1}^d$ is ONB in $\mathbb{R}^d$, with 
$$
\nabla^Tf(x)=P_T(\nabla f(x)),\quad f\in C_c^1(\mathcal{D}),
$$
and $P_T$ is the orthogonal projection of $\mathbb{R}^d$ onto $T$.
\end{definition}

For any $T\in G(d,d-1)$ with $p\in P$ the unit normal of $T$, we have that
\begin{align*}
\text{div}_T\vec{Y}&=\sum_{i=1}^d\nabla^T_iY^i=\sum_{i=1}^de_i\cdot(P_T(\nabla Y^i))\\
&=\sum_{i=1}^de_i\cdot\left(\nabla Y^i-(\nabla Y^i\cdot p)p\right)\\
&=\sum_{i=1}^d\left(\partial_{x_i}Y^i-\sum_{j=1}^d\partial_{x_j}Y^ip_jp_i\right)\\
&=\nabla\vec{Y}:\left(\text{I}-p\otimes p\right).
\end{align*}
 We simply denote $P\equiv G(d,d-1)$. Hence the first variation formula becomes
\begin{equation}\label{2.3}
\langle\delta V,\vec{Y}\rangle=\int\int_{\mathcal{D}\times P}\nabla\vec{Y}(x):\left(\text{I}-p\otimes p\right)dV(x,p),
\end{equation}

Moreover $V$ is said to have \emph{locally bounded first variation} in $\mathcal{D}$ if for each $U$ compactly embedded in $\mathcal{D}$, i.e. $U$ is open in $\mathcal{D}$ and $\bar{U}$ is compact in $\mathcal{D}$, there exists a constant $c>0$ such that
$$
\langle\delta V,\vec{Y}\rangle\leq c\sup_{U}|\vec{Y}|,\quad\forall\vec{Y}\in C_c^1(U,\mathbb{R}^d).
$$
By the general Riesz representation \cite[Theorem 4.1]{Simon:1983ti}, this is equivalent to that there exists a Radon measure $|\delta V|$ on $\mathcal{D}$ characterized by 
$$
|\delta V|(U)=\sup_{\vec{Y}\in C_c(U;\mathbb{R}^d),|\vec{Y}|\leq 1}|\langle \delta V,\vec{Y}\rangle|<\infty.
$$

\vskip.10in
\textbf{Mean curvature vector}

\begin{definition}\label{d2.5}
Let $V$ be a varifold which has locally bounded first variation in $\mathcal{D}$ such that $|\delta V|$ is absolutely continuous  w.r.t. $\|V\|$. A $\|V\|$-measurable vector-valued function $\vec{H}_{V}$ is called a (generalized) mean curvature vector of $V$, if
\begin{equation}\label{2.4}
-\langle\delta V,\vec{Y}\rangle=\langle \|V\|,\vec{H}_V\cdot\vec{Y}\rangle:=\int_\mathcal{D}\vec{H}_V(x)\cdot\vec{Y}(x)d\|V\|(x).
\end{equation}
\end{definition}

\subsection{Definition of a weak solution to the limit of equation (\ref{1.1})}

Now we recall the following definition of weak solutions to the limit of equation (\ref{1.1}) introduced in \cite[Definition 2.1]{Chen:1996bg}:

\begin{definition}\label{d2.1}
A triple $(E,v,V)$ is called a weak solution to the limit of the Hele-Shaw problem (\ref{1.3}) if the following holds:

(\rmnum{1}) $E=\cup_{t\in[0,T]}(\{t\}\times E_t)$ is a subset of $\mathcal{D}_T$ and $\chi_E\in C([0,T];L^1)\cap L^\infty(0,T;BV)$;

(\rmnum{2}) $v\in L^2(0,T;H^1)$;

(\rmnum{3}) $V=V(t,x,p)$ is Radon measure on $\mathcal{D}_T\times P$ and for almost every $t\in[0,T]$, $V^t:=V(t,\cdot,\cdot)$ is a varifold on $\mathcal{D}$, and there exist  Radon measure $\mu^t$ on $\mathcal{\bar{D}}$, $\mu^t$-measurable functions $c_1^t,\cdots,c_d^t$, and $\mu^t$-measurable $P$-valued functions $p_1^t,\cdots,p_d^t$ such that 
\begin{equation}\label{2.5}
0\leq c_i^t\leq1\;\;(i=1,\cdots,d),\quad \sum_{i=1}^dc_i^t\geq 1,\quad\sum_{i=1}^d p_i^t\otimes p_i^t=\mathrm{I}\;\;\mu^t-a.e.,
\end{equation}
\begin{equation}\label{2.6}
2S|D\chi_{E_t}|(x)dx\leq d\mu^t(x)\;\;\left(S=\int_{-1}^1\sqrt{\frac{F(s)}{2}}ds=\frac{2}{3}\right),
\end{equation}
\begin{equation}\label{2.7}
\int\int_{\mathcal{D}\times P}\psi(x,p)dV^t(x,p)=\sum_{i=1}^d\int_\mathcal{D}c_i^t(x)\psi(x,p_i^t(x))d\mu^t(x)\quad\forall\psi\in C_c(\mathcal{D}\times P);
\end{equation}

(\rmnum{4}) For any $t\in(0,T]$ and for almost every $\tau\in(0,t)$,
\begin{equation}\label{2.8}
\int_0^t\int_\mathcal{D}\left(-2\chi_{E_\tau}\partial_t\psi+\nabla v\cdot\nabla \psi\right)dxd\tau=\int_\mathcal{D}2\chi_{E_0}\psi(0,\cdot)\quad\forall\psi\in C_c^1([0,t)\times\bar{\mathcal{D}}),
\end{equation}
\begin{equation}\label{2.9}
-\langle D\chi_{E_t},v\vec{Y}\rangle:=\langle \chi_{E_t},\mathrm{div}(v\vec{Y})\rangle=\frac{1}{2}\langle\delta V^t,\vec{Y}\rangle\quad\forall \vec{Y}\in C_c^1(\mathcal{D};\mathbb{R}^d),
\end{equation}
\begin{equation}\label{2.10}
\mu^t(\bar{\mathcal{D}})+\int_\tau^t\int_\mathcal{D}|\nabla v|^2\leq\mu^\tau(\bar{\mathcal{D}}).
\end{equation}
\end{definition}
\subsection{Main results for $Q$-Wiener noise}\label{ss2.3}
\begin{theorem}\label{t2.3}
Assume that $\sigma\geq \frac{1}{2}$ and (\ref{1.2}) hold. Let $Q$ satisfy (\ref{2.1}) and (\ref{2.2}). Let $(u^\varepsilon,v^\varepsilon)$ be the solution to (\ref{1.1}). Then there exist a probability space $(\tilde{\Omega},\tilde{\mathcal{F}},\{\tilde{\mathcal{F}}_t\}_{t\in[0,T]},\tilde{\mathbb{P}})$, $(\tilde{u}^\varepsilon,\tilde{v}^\varepsilon)\in C([0,T];L^2)\times L^2(0,T;H^1)$ with  $\tilde{\mathbb{P}}\circ\left(\tilde{u}^\varepsilon,\tilde{v}^\varepsilon\right)^{-1}=\mathbb{P}\circ\left(u^\varepsilon,v^\varepsilon\right)^{-1}$ on $C([0,T];L^2)\times L^2(0,T;H^1)$. There also exists a subsequence $\varepsilon_k$ such that as $\varepsilon_k\searrow 0$ the following holds:

(\rmnum{1}) There exists a measurable set $E\subset\tilde{\Omega}\times\mathcal{D}_T$, such that $\chi_{E}$ is $\{\tilde{\mathcal{F}_t}\}$-adapted in $L^2(\mathcal{D})$ and
 for 
$\tilde{\mathbb{P}}-a.s.\; \omega$
$$
\tilde{u}^{\varepsilon_k}(\omega)\to -1+2\chi_{E(\omega)},\quad a.e. \;\;\text{in}\;\;\mathcal{D}_T\;\;\text{and in}\;\;C^\beta([0,T];L^2)
$$
for any $\beta<\frac{1}{12}$ where $E(\omega):=\{(t,x)\in\mathcal{D}_T:(\omega,t,x)\in E\}$;

(\rmnum{2}) There exists $v$ which is weakly measurable in $L^2(0,T;H^1)$, such that for $\tilde{\mathbb{P}}-a.s.\;\omega$
$$
\tilde{v}^{\varepsilon_k}(\omega)\to v(\omega) \quad\text{weakly in}\;\; L^2(0,T;H^1);
$$

(\rmnum{3}) There exist random variables $\mu\in\mathfrak{M}_R$ and $\{\mu_{ij}\}_{i,j=1}^d\in\mathfrak{M}^{d\times d}$ such that for $\tilde{\mathbb{P}}-a.s.\;\omega$
\begin{equation}
\begin{aligned}
e^{\varepsilon_k(\omega)}(\tilde{u}^{\varepsilon_k(\omega)})dxdt
&\to d\mu(\omega,t,x) \quad\text{weakly in}\;\;\mathfrak{M}_R,\\
\varepsilon_k\partial_{x_i}\tilde{u}^{\varepsilon_k}(\omega)\partial_{x_j}\tilde{u}^{\varepsilon_k}(\omega)dxdt&\to d\mu_{ij}(\omega,t,x)\quad\text{weakly in}\;\;\mathfrak{M},\;\;\forall i,j=1,\cdots,d.
\end{aligned}
\end{equation}

(\rmnum{4}) For $\tilde{\mathbb{P}}-a.s.\;\omega$, there exists  Radon measure $V(\omega)$ on $\mathcal{D}_T\times P$,  and $\mu^t(\omega,x)dt=d\mu(\omega,t,x)$ such that for any $t\in(0,T]$ and $\vec{Y}\in C_c^1(\mathcal{D}_t;\mathbb{R}^d)$
\begin{equation}
\int_0^t\langle\delta V^s,\vec{Y}\rangle ds=\int_0^t\int_\mathcal{D}\nabla\vec{Y}:\left(Id\mu(s,x)-\left(\mu_{ij}(s,x)\right)_{d\times d}\right).
\end{equation}

In particular,  for $\tilde{\mathbb{P}}-a.s.\;\omega$, $\left(E(\omega),v(\omega),V(\omega)\right)$ satisfies all the properties in Definition \ref{d2.1} except (\ref{2.10}).
 If $\sigma>\frac{1}{2}$, (\ref{2.10}) holds, thus $\left(E(\omega),v(\omega),V(\omega)\right)$ is a weak solution in the sense of Definition \ref{d2.1}.

\end{theorem}

\begin{theorem}\label{t2.8}
Let $\sigma\geq\frac{1}{2}$,
with the same notations as in Theorem \ref{t2.3},
and suppose that the assumptions in Theorem \ref{t2.3} hold. Then in radially symmetric case, that is
$\mathcal{D}=B_1$, where $B_1$ is the unit ball in $\mathbb{R}^d$ and that $u_0^\varepsilon$ is radially symmetric,  we have that
$$
d\mu=2S|D\chi_{E_t}|dxdt\;\;\text{as Radon measure on}\;\;\mathcal{D}_T.
$$
In particular, for $a.e.t\in[0,T]$, $V^t$ is a $(d-1)$-rectifiable varifold (see \cite[Section 11, Section 38]{Simon:1983ti} for definition), i.e.
$$
dV(t,x,p)=2S|D\chi_{E_t}|dxdt\delta_{\vec{\nu}_{E_t}(t,x)}(dp)\;\;\text{as Radon measure on}\;\;\mathcal{D}_T\times P,
$$
where $\vec{\nu}_{E_t}$ is defined in (\ref{normal}).
Then we have that
\begin{equation}\label{2.15.0}
\left\{
\begin{aligned}
\left(d\mu_{ij}\right)_{d\times d}&=\vec{\nu}_{E_t}\otimes\vec{\nu}_{E_t}d\mu\;\;\text{as Radon measure on}\;\;\mathcal{\bar{D}}_T,\\
v(t,x)&=S\vec{\nu}_{E_t}(x)\cdot \vec{H}_{V^t}(x)\;\;\text{on}\;\;\mathrm{supp}(|D\chi_{E_t}|)\;\;\text{for}\;\;a.e.\;t\in[0,T],
\end{aligned}
\right.
\end{equation}
$\vec{H}_{V^t}$ is the mean curvature vector of $V^t$ defined in Definition \ref{d2.5}
and $\delta_{\vec{\nu}}$ is the Dirac measure concentrated at $\vec{\nu}\in P$.
\end{theorem}

\begin{remark}\label{r2.9}
Since $E_t$ is a BV set for $a.e.\;t\in[0,T]$, by \cite[Theorem 14.3]{Simon:1983ti}, we have that in radial symmetric case, for $a.e.\;t\in[0,T]$
$$
\mu^t=2S|D\chi_{E_t}|=2S\mathcal{H}^{d-1}\lfloor \partial^*E_t.
$$
Here $\mathcal{H}^{d-1}\lfloor \partial^*E_t$ is the $(d-1)$-dimensional Hausdorff measure on $\partial^*E_t$. $\partial^*E_t$ is the so-called reduced boundary of $E_t$ (see  \cite[Section 14]{Simon:1983ti} for details). In our case
$$
\partial^*E_t=\{x\in\mathcal{D}:|\vec{\nu}_{E_t}(x)|=1\}=\mathrm{supp}(|D\chi_{E_t}|).
$$

Since $V^t$ is a rectifiable varifold (see \cite[Chapter 8]{Simon:1983ti} for details), we also have
$$
\vec{H}_{V^t}=\vec{H}_{\partial^*{E_t}},
$$
where $\vec{H}_{\partial^*{E_t}}$ is the generalized mean curvature vector of $\partial^*{E_t}$ (see \cite[Definition 16.5]{Simon:1983ti}). Moreover, when $E_t$ is a  smooth domain, $\vec{H}_{\partial^*{E_t}}$ is just the classical mean curvature vector of $\partial E_t$ and $\vec{\nu}_{E_t}$ is the inward normal vector of $\partial E_t$. Thus the last equation in (\ref{2.15.0}) gives a weak formula of the third equation in (\ref{1.3}).
\end{remark}

\subsection{Remarks on the definition of weak solutions}\label{ss2.4}

Suppose that $(E,v,V)$ is a weak solution of Definition \ref{d2.1}. In the following, we show how Definition \ref{d2.1} is connected with (\ref{1.3}). This has been obtained in \cite[Subsection 2.4]{Chen:1996bg}. We give more details for complete results.

Observe that in distribution sense, $\partial_t\chi_E$ is defined for any $\psi\in C_c^1([0,t)\times\bar{\mathcal{D}})$
$$
\int_0^t\int_\mathcal{D}(\partial_t\chi_E)\psi=\int_0^t\int_\mathcal{D}\partial_t(\chi_E\psi)-\int_0^t\int_\mathcal{D}\chi_E\partial_t\psi=-\int_\mathcal{D}\chi_{E_0}\psi(0,x)dx-\int_0^t\int_\mathcal{D}\chi_E\partial_t\psi dxds,
$$
Thus (\ref{2.8}) implies that in distribution sense
$$
2\partial_t \chi_E=\Delta v,\quad\text{in}\;\;[0,T]\times\mathcal{D}.
$$

Since $v\in L^2(0,T;H^1)$, $\Delta v$ and $\frac{\partial v}{\partial n}$ are ill-defined in (\ref{1.3}). They have to be understood in distribution sense. We suppose that $(v,\Gamma)$ is smooth enough such that $\Delta v$ and $\frac{\partial v}{\partial n}$ are well-defined.

Suppose that $\bar{E}\subset \mathcal{D}$.
Denote $\Gamma_t:=\partial E_t$ and let $\mathcal{D}^+=E_t^o\cap \mathcal{D}$ be the interior of $E_t$ in $\mathcal{D}$ and $\mathcal{D}^-=\mathcal{D}\setminus\bar{E_t}$. 

\textbf{For the first equation in (\ref{1.3}):} For any $x\in\mathcal{D}\setminus\Gamma$, $\Delta v(x)=0$ since $\chi_{E}(x)$ is a constant in time. More precisely,
let $\psi\in C_c^1([0,t)\times\bar{\mathcal{D}})$ and $\mathrm{supp}\psi(s,\cdot)\subset\mathcal{D}\setminus\Gamma_s$ for any $s\in[0,t)$, we have that
$$
\int_0^t\int_\mathcal{D}\chi_{E_t}\partial_t\psi dxds=\int_0^t\int_\mathcal{D}\partial_t\psi dxds
=-\int_\mathcal{D}\psi(0,\cdot)dx=-\int_\mathcal{D}\chi_{E_0}\psi(0,\cdot)dx.
$$
Then (\ref{2.8}) implies that
$$
\int_0^t\int_\mathcal{D}\nabla v\cdot\nabla\psi dxds=0,
$$
which is the weak formula of the first equation in (\ref{1.3}).

\textbf{For the second equation in (\ref{1.3}):} Since $\bar{E}\subset \mathcal{D}$, $\partial \mathcal{D}^-=\partial\mathcal{D}\cup\Gamma$.
For any $\psi\in C_c^1([0,t)\times\bar{\mathcal{D}^-})$ and $\mathrm{supp}\psi(s,\cdot)\subset\mathcal{D}^-$ for any $s\in[0,t)$,
\begin{equation}\label{2.14}
\begin{aligned}
\int_0^t\int_{\partial\mathcal{D}}\frac{\partial v}{\partial n}\psi d\mathcal{H}^{d-1}ds=&\int_0^t\int_\mathcal{D^-}\text{div}(\nabla v\psi)dxds\\
=&\int_0^t\int_{\mathcal{D}^-}\nabla v\cdot\nabla\psi dxds+\int_0^t\int_{\mathcal{D}^-}\Delta v \psi dxds\\
=&\int_0^t\int_\mathcal{D}\nabla v\cdot\nabla\psi dxds+2\int_0^t\int_\mathcal{D}(\partial_t\chi_E)\psi dxds\\
=&\int_0^t\int_\mathcal{D}\nabla v\cdot\nabla\psi dxds-2\int_\mathcal{D}\chi_{E_0}\psi(0,x)dx-2\int_0^t\int_\mathcal{D}\chi_E\partial_t\psi dxds\\
=&0,
\end{aligned}
\end{equation}
where we used (\ref{2.8}) in the last equality. 
Thus we obtain in distribution sense the following holds.
$$
\frac{\partial v}{\partial n}=0,\quad\text{on}\;\;[0,T]\times\partial\mathcal{D}.
$$

\textbf{For the last equation in (\ref{1.3}):} For any $\psi\in C_c^1(\bar{\mathcal{D}}_t)$
\begin{equation}\label{2.15}
\begin{aligned}
\int_0^t\int_\mathcal{D}\partial_t\chi_{E_t}\psi d\mathcal{H}^dds&=-\frac{1}{2}\int_0^t\int_\mathcal{D}\nabla v\nabla \psi d\mathcal{H}^dds\\
&=-\frac{1}{2}\int_0^t\int_\mathcal{D^+}\nabla v\nabla \psi d\mathcal{H}^dds-\frac{1}{2}\int_0^t\int_\mathcal{D^-}\nabla v\nabla \psi d\mathcal{H}^dds\\
&=\frac{1}{2}\int_0^t\int_\mathcal{D^+}\text{div}(\nabla v\psi)d\mathcal{H}^dds+\frac{1}{2}\int_0^t\int_\mathcal{D^-}\text{div}(\nabla v\psi)d\mathcal{H}^dds\\
&=\frac{1}{2}\int_0^t\int_{\Gamma_t}(\partial_n v^+-\partial_n v^-)\psi d\mathcal{H}^{d-1}ds,
\end{aligned}
\end{equation}
which yields that in distribution sense
$$
\mathcal{V}=\frac{1}{2}(\partial_n v^+-\partial_n v^-).
$$

Therefore we know that (\ref{2.8}) is a weak formulation of all the equations in (\ref{1.3}) except the third equation. 

\textbf{For the third equation in (\ref{1.3}):}
following the argument in \cite[Subsection 2.4]{Chen:1996bg},
we can only prove the third equation in weak sense in the radial symmetric case as in Theorem \ref{t2.8} and Remark \ref{r2.9}. 

In general case, it was shown  in \cite[Theorem 3.1, Theorem 3.2]{Roger:2008du}, under the assumption that for $a.e.\;t\in[0,T]$, $v^\varepsilon(t,\cdot)\to v(t,\cdot)$ weakly in $W^{1,p}$ for $p>d$, the authors proved that
\begin{equation}
v(t,x)=S\vec{H}_{\partial^*E_t}\cdot\vec{\nu}_{E_t},\;\;\mathcal{H}^{d-1}-a.e.\;x\in \partial^*E_t.
\end{equation}
But the assumption that $v^\varepsilon\to v$ weakly in $W^{1,p}$ for $p>d$ has not been obtained until now since we can only obtain the convergence in $H^1=W^{1,2}$.

In fact, identifying the value of $v$ on the interface $\Gamma_t$ is the main task of varifold approach to study the sharp interface limit of both Cahn-Hilliard equation and Allen-Cahn equation (cf. \cite{Hutchinson:2000vy,Tonegawa:2002dl,Tonegawa:2005ut,Roger:2006hz,Le:2008dt,Roger:2008du}). Until now, a fully rigorous proof for the (deterministic) Cahn-Hilliard equation is still missing.

\section{Convergence}\label{s3}

\subsection{Lyapunov functional $\mathcal{E}^\varepsilon$ and basic estimates}

In the deterministic case, where no forcing terms are present, the function $\mathcal{E}^\varepsilon$ defined in (\ref{3.1}) decreases in time. In stochastic case, the authors in \cite{DaPrato:1996kk} showed a similar property when $\varepsilon=1$ and (\ref{2.1})
is satisfied. Using the same trick we can prove  a similar result.

\begin{lemma}\label{l3.1}
There exists a constant which only depends on $T$ and $0<\varepsilon_0<1$ such that for any $\varepsilon\in(0,\varepsilon_0]$ and any $p\geq 1$,
\begin{equation}\label{3.2}
\mathbb{E}\sup_{t\in[0,T]}\mathcal{E}^\varepsilon(t)^p\leq C_T(\varepsilon^{2\sigma-1}+\mathcal{E}_0)^p,
\end{equation}
and
\begin{equation}\label{3.3}
\mathbb{E}\left(\int_0^T\|\nabla v^\varepsilon\|_{L^2}^2dt\right)^p \leq C_T(\varepsilon^{2\sigma-1}+\mathcal{E}_0)^p.
\end{equation}
\end{lemma}
\proof
We will not give all the details of the proof since it is the same as \cite[Section 2.3]{DaPrato:1996kk}, we only calculate the order of $\varepsilon$ for every term in the following.

Applying It\^o's formula on $\mathcal{E}^\varepsilon$, we have that
\begin{equation}\label{3.4}
\begin{aligned}
d\mathcal{E}^\varepsilon(u^\varepsilon)&=\langle D\mathcal{E}^\varepsilon(u^\varepsilon),du^\varepsilon\rangle+\frac{\varepsilon^{2\sigma}}{2}\text{Tr}(Q D^2\mathcal{E}^\varepsilon(u^\varepsilon))dt\\
&=-\langle\nabla v^\varepsilon,\nabla v^\varepsilon\rangle dt+\frac{\varepsilon^{2\sigma+1}}{2}\text{Tr}(-\Delta Q)dt+\frac{\varepsilon^{2\sigma-1}}{2}\text{Tr}(f'(u^\varepsilon)Q)dt+\varepsilon^\sigma\langle v^\varepsilon,dW_t\rangle.
\end{aligned}
\end{equation}
By using the same trick as in \cite[Section 2.3]{DaPrato:1996kk} we have that 
$$
\text{Tr}(f'(u^\varepsilon)Q)\lesssim 1+\varepsilon\mathcal{E}^\varepsilon(u^\varepsilon).
$$
Hence we deduce from (\ref{3.4}) that for any $p\geq 1$,
\begin{align*}
\mathbb{E}\left(\sup_{t\in[0,T]}\mathcal{E}^\varepsilon(t)+\int_0^T\|\nabla v^\varepsilon\|_{L^2}^2ds\right)^p
&\lesssim\mathbb{E}\left(\mathcal{E}_0+\varepsilon^{2\sigma-1}+\varepsilon^{2\sigma+1}+\varepsilon^{2\sigma}\sup_{t\in[0,T]}\mathcal{E}^\varepsilon(t)+\varepsilon^{\sigma}\sup_{t\in[0,T]}|M^\varepsilon(t)|\right)^p
\end{align*}
where $M^\varepsilon(t):=\int_0^t\langle v^\varepsilon,dW_s\rangle$. 
Let $\varepsilon$ be small enough, we have that
$$
\mathbb{E}\left(\sup_{t\in[0,T]}\mathcal{E}^\varepsilon(t)\right)^p+\mathbb{E}\left(\int_0^T\|\nabla v^\varepsilon\|_{L^2}^2ds\right)^p\lesssim \varepsilon^{p(2\sigma-1)}+\mathcal{E}_0^p+\mathbb{E}\sup_{t\in[0,T]}|M^\varepsilon(t)|^p,
$$

By Burkholder-Davis-Gundy's inequality
\begin{align*}
\mathbb{E}\sup_{t\in[0,T]}|M^\varepsilon(t)|^p
&\lesssim \mathbb{E}\left(\langle M^\varepsilon\rangle_T\right)^{\frac{p}{2}}=\mathbb{E}\left(\int_0^T\|\sqrt{Q}v^\varepsilon(t)\|_{L^2}^2dt\right)^\frac{p}{2}\\
&\lesssim \mathbb{E}\left(\int_0^T\|\nabla v^\varepsilon\|_{L^2}^2dt\right)^{\frac{p}{2}}.
\end{align*}
Then by Young's inequality, for any $\kappa>0$, there exists a constant $C_1\equiv C_1(T)$ such that
$$
\mathbb{E}\sup_{t\in[0,T]}|M^\varepsilon(t)|^p\leq C_1+\kappa\mathbb{E}\left(\int_0^T\|\nabla v^\varepsilon\|_{L^2}^2ds\right)^p.
$$
Thus for a small enough $\kappa>0$, there exists a constant $C_T>0$ such that
$$
\mathbb{E}\sup_{t\in[0,T]}\mathcal{E}^\varepsilon(t)^p+(1-\kappa)\mathbb{E}\left(\int_0^T\|\nabla v^\varepsilon\|_{L^2}^2ds\right)^p\leq C_T(\varepsilon^{p(2\sigma-1)}+\mathcal{E}_0^p),
$$
which implies our results.
$\hfill\Box$
\vskip.10in

\begin{corollary}\label{c3.2}
There exists a constant $C_T>0$, such that for any $p\geq 1$
\begin{equation}\label{3.6}
\mathbb{E}\sup_{t\in[0,T]}\left(\int_\mathcal{D}F(u^\varepsilon(t,x))dx\right)^p\leq C_T\varepsilon^p(\mathcal{E}_0^p+\varepsilon^{p(2\sigma-1)})
\end{equation}
and
\begin{equation}\label{3.7}
\mathbb{E}\sup_{t\in[0,T]}\|u^\varepsilon(t)\|_{L^4}^{4p}\leq C_T(1+\mathcal{E}_0^p+\varepsilon^{p(2\sigma-1)}).
\end{equation}
\end{corollary}

In the rest of this section, we always assume $\sigma\geq \frac{1}{2}$.

\subsection{Estimates for $\{u^\varepsilon\}$}\label{ss3.3}
We introduce a function $g^\varepsilon(t,x)$ defined by
\begin{equation}\label{3.8}
g^\varepsilon(t,x):=G(u^\varepsilon(t,x)),
\end{equation}
where
$$
G(u):=\int_{-1}^u\sqrt{2F(x)}dx,\quad\forall u\in\mathbb{R}.
$$
Observe that
\begin{equation}\label{3.9}
\int_\mathcal{D}|\nabla g^\varepsilon(t,\cdot)|=\int_\mathcal{D}\sqrt{2F(u^\varepsilon)}|\nabla u^\varepsilon|dx\leq\int_\mathcal{D}e^\varepsilon(u^\varepsilon)(t)dx=\mathcal{E}^\varepsilon(t),
\end{equation}
and there are positive constants $c_1$, $c_2$ such that
\begin{equation}\label{3.10}
c_1|u_1-u_2|^2\leq|G(u_1)-G(u_2)|\leq c_2|u_1-u_2|(1+|u_1|+|u_2|),\quad\forall u_1, u_2\in\mathbb{R}.
\end{equation}
\begin{lemma}\label{l3.3}
There exists constant $C_T>0$ which only depends on $T$, such that for any $\beta\in(0,\frac{1}{12})$,
$$
\mathbb{E}\left(\|g^\varepsilon\|_{L^\infty(0,T;W^{1,1})}+\|g^\varepsilon\|_{C^\beta([0,T];L^1)}+\|u^\varepsilon\|_{C^\beta([0,T];L^2)}\right)\leq C_T
$$
\end{lemma}
\proof
Similarly to the proof of \cite[Lemma 3.2]{Chen:1996bg}, let $\rho$ be any fixed mollifier satisfying
$$
\rho\in C^\infty(\mathbb{R}^d),\quad 0\leq\rho\leq 1,\quad\text{supp}\rho\subset B_1(0),\quad\int_{\mathbb{R}^d}\rho=1,
$$
where $B_1$ is the unit ball in $\mathbb{R}^d$ centered at $0$. For any small $\eta>0$, we define
$$
u_\eta^\varepsilon(t,x)=\int_{B_1}\rho(y)u^\varepsilon(t,x-\eta y)dy.
$$
Here we assume that $u^\varepsilon$ is extended to $\{x\in\mathbb{R}^d:d(x,\mathcal{D})\leq\eta_0\}$ by
$$
u^\varepsilon(t,y+\eta n(y))=u^\varepsilon(t,y-\eta n(y)),\quad y\in\partial\mathcal{D},\eta\in[0,\eta_0],
$$
where $\eta_0$ is a small positive number and $n(y)$ is the unit outward normal vector to $\partial\mathcal{D}$ at $y\in\partial\mathcal{D}$.

Then by (\ref{3.7}), we have that for any $p>1$, $\eta\in(0,\eta_0)$,
\begin{equation}\label{3.11}
\mathbb{E}\sup_{t\in[0,T]}\|\nabla u^\varepsilon_\eta(t)\|_{L^2}^p\lesssim \eta^{-p}\mathbb{E}\sup_{t\in[0,T]}\|u^\varepsilon(t)\|_{L^2}^p\lesssim\eta^{-p},
\end{equation}
and
\begin{equation}\label{3.12}
\begin{aligned}
\mathbb{E}\sup_{t\in[0,T]}\left(\int_\mathcal{D}|u_\eta^\varepsilon-u^\varepsilon|^2dx\right)^p&\leq\mathbb{E}\left(\sup_{t\in[0,T]}\int_\mathcal{D}\int_{B_1}\rho(y)|u^\varepsilon(t,x-\eta y)-u^\varepsilon(t,x)|^2dydx\right)^p\\
&\lesssim\mathbb{E}\left(\sup_{t\in[0,T]}\int_\mathcal{D}\int_{B_1}\rho(y)|g^\varepsilon(t,x-\eta y)-g^\varepsilon(t,x)|dydx\right)^p\\
&\lesssim\eta^p\mathbb{E}\sup_{t\in[0,T]}\|\nabla g^\varepsilon(t)\|_{L^1}^p\\
&\leq\eta^p\mathbb{E}\sup_{t\in[0,T]}\mathcal{E}^\varepsilon(t)^p\lesssim \eta^p,
\end{aligned}
\end{equation}
where we use (\ref{3.10}) in the second inequality and (\ref{3.9}), (\ref{3.2}) in the last line.

For any $0\leq\tau<t\leq T$, by using $u^\varepsilon(t)-u^\varepsilon(\tau)=\int_\tau^t\Delta v^\varepsilon(s)ds+\varepsilon^\sigma(W_t-W_\tau)$ (in weak sense), we have that
\begin{equation}\label{3.13}
\begin{aligned}
&\quad\mathbb{E}\left(\int_\mathcal{D}|\left( u^\varepsilon_\eta(t,x)-u^\varepsilon_\eta(\tau,x)\right)\left(u^\varepsilon(t,x)-u^\varepsilon(\tau,x)\right)|dx\right)^p\\
\leq&\mathbb{E}\left(\int_\tau^t\int_\mathcal{D}|\nabla v^\varepsilon(s,x)\left(\nabla u^\varepsilon_\eta(t,x)-\nabla u^\varepsilon_\eta(\tau,x)\right)|dxds\right)^p\\
&+\varepsilon^{p\sigma}\mathbb{E}\left(\int_\mathcal{D}|\left(u_\eta^\varepsilon(t,x)-u_\eta^\varepsilon(\tau,x)\right)\left(W_t-W_\tau\right)|dx\right)^p\\
\lesssim&\mathbb{E}\left(\int_\tau^t\int_\mathcal{D}|\nabla v^\varepsilon|^2\right)^{\frac{p}{2}}(t-\tau)^{\frac{p}{2}}\left(\mathbb{E}\sup_{s\in[0,T]}\|\nabla u_\eta^\varepsilon(s)\|_{L^2}^p\right)\\
&+\varepsilon^{p\sigma}\mathbb{E}\sup_{s\in[0,T]}\|u_\eta^\varepsilon(s)\|_{L^2}^p\left(\mathbb{E}\|W_t-W_\tau\|_{L^2}^{2p}\right)^{\frac{1}{2}}\\
\lesssim&(t-\tau)^{\frac{p}{2}}\eta^{-p}+(t-\tau)^{\frac{p}{2}}\varepsilon^{p\sigma}\\
\lesssim&\eta^{-p}(t-\tau)^\frac{p}{2},
\end{aligned}
\end{equation}
where in the third inequality we use (\ref{3.3}), (\ref{3.7}), (\ref{3.11}) and the fact that
$$
\mathbb{E}\|W_t-W_\tau\|_{L^2}^{2p}\lesssim|t-\tau|^p.
$$
Then we have that
\begin{equation*}
\begin{aligned}
\mathbb{E}\|u^\varepsilon(t)-u^\varepsilon(\tau)\|_{L^2}^{2p}\lesssim&\mathbb{E}\left(\int_\mathcal{D}|\left( u^\varepsilon_\eta(t,x)-u^\varepsilon_\eta(\tau,x)\right)\left(u^\varepsilon(t,x)-u^\varepsilon(\tau,x)\right)|dx\right)^p\\
&+\mathbb{E}\left(\int_\mathcal{D}|\left( u^\varepsilon(t,x)-u^\varepsilon_\eta(t,x)\right)\left(u^\varepsilon(t,x)-u^\varepsilon(\tau,x)\right)|dx\right)^p\\
&+\mathbb{E}\left(\int_\mathcal{D}|\left( u^\varepsilon(\tau,x)-u^\varepsilon_\eta(\tau,x)\right)\left(u^\varepsilon(t,x)-u^\varepsilon(\tau,x)\right)|dx\right)^p\\
\lesssim&\eta^{-p}(t-\tau)^\frac{p}{2}+\left(\mathbb{E}\left(\sup_{t\in[0,T]}\|u_\eta^\varepsilon(t)-u^\varepsilon(t)\|_{L^2}\right)^{2p}\right)^{\frac{1}{2}}\left(\mathbb{E}\sup_{t\in[0,T]}\|u^\varepsilon\|_{L^2}^{2p}\right)^{\frac{1}{2}}\\
\lesssim&\eta^{-p}(t-\tau)^\frac{p}{2}+\eta^{\frac{p}{2}},
\end{aligned}
\end{equation*}
where we use (\ref{3.13}) in the second inequality and (\ref{3.7}), (\ref{3.12}) in the last inequality. If we take $\eta=\eta_0\wedge(t-\tau)^{\frac{1}{3}}$, we have that
\begin{equation}\label{3.14}
\mathbb{E}\|u^\varepsilon(t)-u^\varepsilon(\tau)\|_{L^2}^{2p}\lesssim\eta^{-p}(t-\tau)^\frac{p}{2}+\eta^{\frac{p}{2}}\leq(t-\tau)^{\frac{p}{6}}.
\end{equation}
Moreover, using (\ref{3.10}) we have that
\begin{equation}\label{3.15}
\begin{aligned}
\mathbb{E}\|g^{\varepsilon}(t)-g^\varepsilon(\tau)\|_{L^1}^p&\lesssim\mathbb{E}\left(\int_\mathcal{D}|u^\varepsilon(t,x)-u^\varepsilon(\tau,x)|\left(1+|u^\varepsilon(t,x)|+|u^\varepsilon(\tau,x)|\right)dx\right)^p\\
&\lesssim\mathbb{E}\|u^\varepsilon(t)-u^\varepsilon(\tau)\|_{L^2}^p\left(1+\mathbb{E}\sup_{t\in[0,T]}\|u^\varepsilon\|_{L^2}^p\right)\\
&\lesssim(t-\tau)^\frac{p}{12},
\end{aligned}
\end{equation}
where we use (\ref{3.7}) and (\ref{3.14}) in the last inequality.

Finally by Kolmogorov's criteria (see e.g. \cite[Theorem 3.3]{DaPrato:2014jd}), for any $0<\beta<\frac{1}{12}$,
$$
\mathbb{E}\left(\|g^\varepsilon\|_{C^\beta([0,T];L^1)}+\|u^\varepsilon\|_{C^\beta([0,T];L^2)}\right)\lesssim 1.
$$

Moreover by (\ref{3.9})
$$
\mathbb{E}\sup_{t\in[0,T]}\|\nabla g^\varepsilon(t)\|_{L^1}\lesssim 1.
$$
Thus
$$
\mathbb{E}\|g\|_{L^\infty(0,T;W^{1,1})}\lesssim\mathbb{E}\sup_{t\in[0,T]}\|\nabla g^\varepsilon(t)\|_{L^1}+\mathbb{E}\sup_{t\in[0,T]}\|g^\varepsilon(t)\|_{L^1}\lesssim 1
$$

$\hfill\Box$
\vskip.10in

\subsection{Estimates for $\{v^\varepsilon\}$}\label{ss3.4}
We want to obtain the estimate of $v^\varepsilon$ in the space $H^1$. By (\ref{3.3}) and Poincar\'e-Wirtinger inequality, it is enough to estimate $\bar{v^\varepsilon}:=\frac{1}{|\mathcal{D}|}\int_\mathcal{D}v^\varepsilon(x)dx$.
\begin{lemma}\label{l3.4}
For any $\delta>0$, there exists a constant $C\equiv C(\delta,T)>0$, such that
$$
\mathbb{P}\left(\int_0^T\|v^\varepsilon(t)\|_{H^1}^2dt\leq C\right)\geq 1-\delta.
$$
\end{lemma}
\proof
For any $R>0$, set
$$A_R:=\left\{\omega\in\Omega:\|u^\varepsilon(\omega)\|_{C([0,T];L^2)}+\sup_{t\in[0,T]}\mathcal{E}^\varepsilon(t)(\omega)^p+\int_0^T\|\nabla v^\varepsilon(\omega,t)\|_{L^2}^2dt\leq R\right\}.$$
By the same argument as in \cite[Lemma 3.4]{Chen:1996bg} and
using an integration by parts formula, we have that
$$
\bar{v}^\varepsilon=\frac{\int_\mathcal{D}\left(D^2\psi:\left(e(u^\varepsilon)\mathrm{I}-\varepsilon\nabla u^\varepsilon\otimes\nabla u^\varepsilon\right)-u^\varepsilon\nabla\psi\cdot\nabla v^\varepsilon-u^\varepsilon\Delta\psi(v^\varepsilon-\bar{v}^\varepsilon)\right)}{\int_\mathcal{D}\Delta\psi u^\varepsilon}，
$$
where $D^2\psi$ is the Hessen matrix of $\psi$, $\psi$ is the unique solution to 
\begin{equation*}
\left\{
\begin{aligned}
-\Delta\psi&=u_\eta^\varepsilon-\bar{u}^\varepsilon_\eta\;\;\text{in}\;\;\mathcal{D},\\
\frac{\partial\psi}{\partial n}&=0\;\;\text{on}\;\;\partial\mathcal{D}.
\end{aligned}
\right.
\end{equation*}
Here $u_\eta^\varepsilon$ is defined in the same way as in the proof of Lemma \ref{l3.3}.

Then for a fixed $\omega\in A_R$, all the estimates in the proof of \cite[Lemma 3.4]{Chen:1996bg} hold. Thus we have that there exists a constant $C_R$ such that for any $\omega\in A_R$, $t\in[0,T]$ and any $\eta,\varepsilon\in(0,1)$
$$
|\bar{v}^\varepsilon(\omega,t)|\leq C_R\frac{\eta^{-1}(1+\varepsilon^{1/2}\eta^{-d/2})(\mathcal{E}^\varepsilon(t)(\omega)+\|\nabla v^\varepsilon(\omega,t)\|_{L^2(\mathcal{D})})}{1-m_0^2-\sqrt{\varepsilon}-\sqrt{\eta}},
$$
where $m_0=\bar{u^\varepsilon_0}\in(-1,1)$ is as in (\ref{1.2}).
Taking $\eta$ small and independent of $\varepsilon$, we obtain that there exists constant $\tilde{C}_R>0$ such that for any $\omega\in A_R$, $t\in[0,T]$,
$$
\int_0^T\bar{v}^\varepsilon(t)^2 dt\leq \tilde{C}_R.
$$
Hence we have
$$
A_R\subset\left\{\int_0^T\bar{v}^\varepsilon(t)^2 dt\leq \tilde{C}_R\right\}.
$$
Moreover, by Poincar\'e-Wirtinger inequality
$$
\|v^\varepsilon\|_{H^1}\lesssim|\bar{v^\varepsilon}|+\|\nabla v^\varepsilon\|_{L^2},
$$
thus for any $R>0$ there exists a constant $C_R>0$, such that
$$
\mathbb{P}\left(\int_0^T\|v^\varepsilon(t)\|_{H^1}^2dt\leq C_R\right)\geq \mathbb{P}\left(\int_0^T\bar{v}^\varepsilon(t)^2 dt\leq \tilde{C}_R,\|\nabla v^\varepsilon\|_{L^2(\mathcal{D_T})}^2\leq R\right)\geq \mathbb{P}(A_R).
$$

By Lemma \ref{l3.1} and Lemma \ref{l3.3}, using Cheybeshev's inequality, we have that for any $\delta>0$, there exists a constant $R\equiv R(\delta)>0$, such that
$$
\mathbb{P}\left(A_R\right)\geq 1-\delta.
$$
Then we obtain the assertion of the lemma.
$\hfill\Box$
\vskip.10in

\subsection{Tightness}
For any $\beta<\frac{1}{12}$, we denote
\begin{equation}
\mathcal{X}^1:=\mathbb{R}\times L_{w^\ast}^\infty(0,T)\times C^\beta([0,T];L_w^2)\times C^\beta([0,T];L^1)\times L^2_w(0,T;H^1),
\end{equation}
where $L^2_w(0,T;H^1)$ is the space $L^2(0,T;H^1)$ equipped with the weak topology, $L_w^2$ is the space $L^2$ equipped with the weak topology and $ L_{w^\ast}^\infty(0,T)$ is the space $L^\infty(0,T)$ equipped with the weak-* topology. We also denote
\begin{equation}
\mathcal{X}^2:=\mathfrak{M}^{d\times d}\times\mathfrak{M}_R,
\end{equation}
where $\mathfrak{M}$ is the space of all finite signed measure on $\mathcal{D}_T$ and $\mathfrak{M}_R\subset\mathfrak{M}$ is the space of all Radon measure on $\mathcal{D}_T$. $\mathfrak{M}_R$ and $\mathfrak{M}$ are equipped with the total variation norm $\|\cdot\|_{TV}$ and weak topology, respectively. Here an element in $\mathfrak{M}^{d\times d}$ is a $d\times d$ $\mathfrak{M}$-valued matrix $\{\mu_{ij}\}_{i,j=1}^d$ where $\mu_{ij}\in\mathfrak{M}$.

Let $\hat{\mathbb{P}}^\varepsilon$ be the probability measure on $\mathcal{X}^1\times\mathcal{X}^2$ defined by
\begin{equation}
\hat{\mathbb{P}}^\varepsilon:=\mathbb{P}\circ\left(\varepsilon^{-1}\sup_{t\in[0,T]}\|F(u^\varepsilon)\|_{L^1},\mathcal{E}^\varepsilon(u^\varepsilon),u^\varepsilon,G(u^\varepsilon),v^\varepsilon,e^\varepsilon(u^\varepsilon)dxdt,\{\varepsilon\partial_{x_i}u^\varepsilon\partial_{x_j}u^\varepsilon dxdt\}_{ij}\right)^{-1}.
\end{equation}
In the following we will prove that $\{\hat{\mathbb{P}}^\varepsilon\}_{\varepsilon}$ is tight on $\mathcal{X}^1\times\mathcal{X}^2$. This is equivalent to prove the tightness of every component.

For $\sup_{t\in[0,T]}\|F(u^\varepsilon)\|_{L^1}$, by (\ref{3.6}) and Chebyshev's inequality, we know that
$$
\mathbb{E}\varepsilon^{-1}\sup_{t\in[0,T]}\int_\mathcal{D}F(u^\varepsilon)dx\lesssim 1.
$$
Then we have that for any $\delta>0$, there exists a constant $K_1>0$ such that
$$
\mathbb{P}\left(\varepsilon^{-1}\sup_{t\in[0,T]}\|F(u^\varepsilon)\|_{L^1}\leq K_1\right)\geq 1-\delta.
$$

For $\mathcal{E}^\varepsilon$, by (\ref{3.2}) and Chebyshev's inequality, we have that for any $\delta>0$, there exists a constant $K_2>0$ such that
$$
\mathbb{P}\left(\sup_t \mathcal{E}_t^\varepsilon\leq K_2\right)\geq 1-\delta.
$$
By the Banach-Alaoglu theorem, any bounded set in $L^\infty(0,T)$ is precompact in $L^\infty_{w^\ast}(0,T)$, thus $\mathbb{P}\circ\left(\mathcal{E}^\varepsilon(u^\varepsilon)\right)^{-1}$ is tight on $L^\infty_{w^\ast}(0,T)$.

For $u^\varepsilon$, by the Banach-Alaoglu theorem, any bounded set in $L^2$ is precompact in $L^2_w$. Then by a generalized Arzel\`a-Ascoli theorem, any bounded set in $C^\beta([0,T];L^2)$ is precompact in $C^\gamma([0,T];L^2_w)$ for any $0<\gamma<\beta$. Hence we obtain the tightness of $\mathbb{P}\circ\left(u^\varepsilon\right)^{-1}$ on $C^\gamma([0,T];L^2_w)$ by using Chebyshev' inequality and Lemma \ref{l3.3}.

For $G(u^\varepsilon)$, by Lemma \ref{l3.3} we have that for any $\delta>0$, there exists a constant $K_3>0$ such that
$$
\mathbb{P}\left(\|G(u^\varepsilon)\|_{L^\infty(0,T;W^{1,1})}+\|G(u^\varepsilon)\|_{C^\beta([0,T];L^1)}\leq K_3\right)\geq 1-\delta.
$$
Since $W^{1,1}$ is compactly embedded into $L^q$ for any $q\in[1,\frac{d}{d-1}]$, then by a generalized Arzel\`a-Ascoli theorem for any $0<\gamma<\beta$, the set
$$
\{g\in C^\gamma([0,T];L^1):\|g\|_{L^\infty(0,T;W^{1,1})}+\|g\|_{C^\beta([0,T];L^1)}\leq K\}
$$
is compact in $C^\gamma([0,T];L^1)$, which implies the tightness of $\mathbb{P}\circ\left(G(u^\varepsilon)\right)^{-1} $ in $C^\gamma([0,T];L^1)$ for any $\gamma<\frac{1}{12}$.

For $v^\varepsilon$, the tightness of $\mathbb{P}\circ\left(v^\varepsilon\right)^{-1} $ in $L^2_w(0,T;H^1)$ is followed by Lemma \ref{l3.4} and the Banach-Alaoglu theorem.

For $\varepsilon\partial_{x_i}u^\varepsilon\partial_{x_j}u^\varepsilon$ and $e^\varepsilon(u^\varepsilon)$, since $L^1(\mathcal{D}_T)$ is embedded into $\mathfrak{M}$. Moreover for any $f\in L^1(\mathcal{D}_T)$, we have that
$$
f(t,x)dxdt=f^+dxdt-f^-dxdt.
$$
Since $\mathcal{D}_T$ is a compact set,  we have that $f^+dxdt, f^-dxdt\in\mathfrak{M}_R$. By \cite[Theorem 4.4]{Simon:1983ti}, 
any bounded set in $\mathfrak{M}_R$ w.r.t. total variation norm is precompact in $\mathfrak{M}_R$ w.r.t weak topology, which implies that any bounded set in $\mathfrak{M}$ w.r.t. total variation norm is precompact in $\mathfrak{M}$ w.r.t weak topology. Thus by
(\ref{3.2}) and 
$$
\|\varepsilon\partial_{x_i}u^\varepsilon\partial_{x_j}u^\varepsilon\|_{L^1(\mathcal{D}_T)}\lesssim\varepsilon\|\nabla u^\varepsilon\|_{L^1(\mathcal{D}_T)}\lesssim\sup_{t\in[0,T]}\mathcal{E}^\varepsilon_t,
$$
$$
\|e^\varepsilon(u^\varepsilon)\|_{L^1(\mathcal{D}_T)}\lesssim\sup_{t\in[0,T]}\mathcal{E}^\varepsilon_t,
$$
we obtain the tightness of $\mathbb{P}\circ\left(e^\varepsilon(u^\varepsilon)dxdt,\{\varepsilon\partial_{x_i}u^\varepsilon\partial_{x_j}u^\varepsilon dxdt\}_{ij}\right)^{-1}$ in $\mathcal{X}^2$.

Hence we proved the tightness of $\{\hat{\mathbb{P}}^\varepsilon\}_\varepsilon$ in $\mathcal{X}^1\times\mathcal{X}^2$. Then by using a Jakubowski’s version of the Skorokhod Theorem in the form given by \cite[Theorem A.1]{Brzezniak:2013bk}, which was proved in \cite{Jak:1998cb}:

\begin{theorem}\label{tsk}
Let $\mathcal{X}$ be a topological space such that there exists a sequence $\{f_n\}_{n\geq 1}$ of continuous functions $f_n : \mathcal{X}\to\mathbb{R}$ that separate points of $\mathcal{X}$. Let us denote by $\mathscr{S}$ the $\sigma$-algebra generated by the maps $\{f_n\}$. Then:

($\mathrm{j1}$) every compact subset of $\mathcal{X}$ is metrizable;

($\mathrm{j2}$) every Borel subset of a $\sigma$-compact set in $\mathcal{X}$ belongs to $\mathscr{S}$ ;

($\mathrm{j3}$) every probability measure supported by a $\sigma$-compact set in $\mathcal{X}$ has a unique
Radon extension to the Borel $\sigma$-algebra on $\mathcal{X}$;

($\mathrm{j4}$) if $(\mu_n)$ is a tight sequence of probability measures on $(\mathcal{X}, \mathscr{S})$, then there
exists a subsequence $(n_k)_{k\geq 1}$, a probability space $( \Omega,\mathcal{F},\mathbb{P})$ with $\mathcal{X}$-valued Borel measurable random variables $X_k$,  $X$ such that $\mu_{n_k}$ is the law of $X_k$ and $X_k$ converge almost surely to $X$. Moreover, the law of $X$ is a Radon measure.
\end{theorem}

We obtain that
\begin{theorem}\label{t3.5} Assume $
\sigma\geq \frac{1}{2}$.
There exist a probability space $(\tilde{\Omega},\tilde{\mathcal{F}},\{\tilde{\mathcal{F}}\}_{t\in[0,T]},\tilde{\mathbb{P}})$ on $\mathcal{X}^1\times\mathcal{X}^2$, a subsequence $\varepsilon_k$ (we still denote it as $\varepsilon$ for simplicity) and 
$$\left\{\left(\varepsilon^{-1}\sup_{t\in[0,T]}\|F(\tilde{u}^\varepsilon)\|_{L^1},\mathcal{E}^\varepsilon(\tilde{u}^\varepsilon),\tilde{u}^\varepsilon,G(\tilde{u}^\varepsilon),\tilde{v}^\varepsilon,e^\varepsilon(\tilde{u}^\varepsilon)dxdt,\{\varepsilon\partial_{x_i}\tilde{u}^\varepsilon\partial_{x_j}\tilde{u}^\varepsilon dxdt\}_{ij}\right)\right\}\subset \mathcal{X}^1\times\mathcal{X}^2$$ 
and 
$$\left(a,\mathcal{E},u,g,v,\mu,\{\mu_{ij}\}_{ij}\right)\in \mathcal{X}^1\times\mathcal{X}^2,$$ 
such that

(\rmnum{1}) $\tilde{P}\circ \left(\varepsilon^{-1}\sup_{t\in[0,T]}\|F(\tilde{u}^\varepsilon)\|_{L^1},\mathcal{E}^\varepsilon(\tilde{u}^\varepsilon),\tilde{u}^\varepsilon,G(\tilde{u}^\varepsilon),\tilde{v}^\varepsilon,e^\varepsilon(\tilde{u}^\varepsilon)dxdt,\{\varepsilon\partial_{x_i}\tilde{u}^\varepsilon\partial_{x_j}\tilde{u}^\varepsilon dxdt\}_{ij}\right)^{-1}=\hat{\mathbb{P}}^\varepsilon$ on $\mathcal{X}^1\times\mathcal{X}^2$,

(\rmnum{2}) $\left(\varepsilon^{-1}\sup_{t\in[0,T]}\|F(\tilde{u}^\varepsilon)\|_{L^1},\mathcal{E}^\varepsilon(\tilde{u}^\varepsilon),\tilde{u}^\varepsilon,G(\tilde{u}^\varepsilon),\tilde{v}^\varepsilon,e^\varepsilon(\tilde{u}^\varepsilon)dxdt,\{\varepsilon\partial_{x_i}\tilde{u}^\varepsilon\partial_{x_j}\tilde{u}^\varepsilon dxdt\}_{ij}\right)$ converges to $\left(0,\mathcal{E},u,g,v,\mu,\{\mu_{ij}\}_{ij}\right)$ in $\mathcal{X}^1\times\mathcal{X}^2$, $\tilde{\mathbb{P}}-a.s$, as $\varepsilon\searrow 0$.

In particular, for $\tilde{\mathbb{P}}-a.s. \omega$, there exists a Borel set $E(\omega)\subset \tilde{\Omega}\times\mathcal{D}_T$, such that as $\varepsilon\searrow 0$

(\rmnum{3}) $u^\varepsilon\to u$ in $C^{\frac{\beta}{2}}([0,T];L^2)$, $g=G(u)=2S\chi_E$ a.e. in $\mathcal{D}_T$ and in $C^\beta([0,T];L^1)$, $u=-1+2\chi_E$ a.e. in $\mathcal{D}_T$ and in $C^\beta([0,T];L^2)$.

Moreover, denote $E=\{(\omega,t,x)\in \Omega\times D_T:(t,x)\in E(\omega)\}$,  $E_t:=\left\{(\omega,x):(\omega,t,x)\in E\right\}$, then $\chi_{E_t}$ is $\{\tilde{\mathcal{F}}\}_{t\in[0,T]}$-adapted in $L^2(\mathcal{D})$ and satisfies the following:

(\rmnum{4}) For all $\beta<\frac{1}{12}$, $\tilde{\mathbb{P}}\left(\chi_{E}\in C^\beta([0,T];L^1)\right)=1$,

(\rmnum{5}) $\tilde{\mathbb{P}}\left(|E_t|=|E_0|=\frac{1+m_0}{2}|\mathcal{D}|,\forall t\in[0,T]\right)=1$,

(\rmnum{6}) $\tilde{\mathbb{P}}\left(\chi_{E}\in L^\infty(0,T;BV)\right)=1$.
\end{theorem}
\proof
Since $\mathcal{X}^1\times\mathcal{X}^2$ is locally convex space and its dual space is separable, by \cite[Theorem 3.4]{Rudin:1973ul}, the condition in Theorem \ref{tsk} holds. Thus
the Skorohod theorem Theorem \ref{tsk} yields the first assertion and the existence of convergence subsequence to 
$$
\left(a,\mathcal{E},u,g,v,\mu,\{\mu_{ij}\}_{ij}\right)\;\; \text{in}\;\; \mathcal{X}^1\times\mathcal{X}^2.
$$

Since $\tilde{\mathbb{P}}\circ(\tilde{u}^\varepsilon,\tilde{v}^\varepsilon)^{-1}=\mathbb{P}\circ(u^\varepsilon,v^\varepsilon)^{-1}$, we have that for any $h\in H^1$,
$$
\varepsilon^{-\sigma}\left(\int_\mathcal{D}(\tilde{u}^\varepsilon(t)-\tilde{u}^\varepsilon(0))hdx+\int_0^t\nabla\tilde{v}^\varepsilon\cdot\nabla hdx\right)
$$
is a Wiener process on $(\tilde{\Omega},\tilde{\mathcal{F}},\tilde{\mathbb{P}})$ with covariance $\|Q^{\frac{1}{2}}h\|_{L^2}^2$. Thus there exists a $Q$-Wiener process $\tilde{W}$ on $L^2$ which is defined on $(\tilde{\Omega},\tilde{\mathcal{F}},\tilde{\mathbb{P}})$. Then we have that for any $h\in H^1$
\begin{equation}\label{3.17.1}
\int_\mathcal{D}(\tilde{u}^\varepsilon(t)-\tilde{u}^\varepsilon(0))hdx+\int_0^t\int_\mathcal{D}\nabla\tilde{v}^\varepsilon\cdot\nabla hdx=\varepsilon^\sigma\int_0^t\langle h,d\tilde{W}_s\rangle.
\end{equation}

Moreover, we denote $\{\tilde{\mathcal{F}}_t\}_{t\in[0,T]}$ be completion under $\tilde{\mathbb{P}}$ of the natural filtration generated by $\{\tilde{W}_t\}_{t\in[0,T]}$, thus $\{\tilde{\mathcal{F}}_t\}_{t\in[0,T]}$ is a normal filtration. By \cite{DaPrato:1996kk}, we know that for any $\varepsilon>0$, $\tilde{u}^\varepsilon$ is the unique solution, thus by Yamada–Watanabe theorem (see e.g. \cite[Theorem E.0.8]{Liu:2015vb}) $\{\tilde{u}^\varepsilon\}_{t}$ is  $\{\tilde{\mathcal{F}}_t\}_{t}$-adapted in $L^2(\mathcal{D})$. Since $\tilde{u}^\varepsilon\to u$ in $C([0,T];L_w^2)$, we know that $u$ is $\{\tilde{\mathcal{F}}_t\}_{t}$-adapted in $L^2(\mathcal{D})$. 

In the rest of this proof, we ignore the notation $\;\tilde{}\;$ if there is no confusion.

By (\ref{3.10}), we know that for any $t,\tau\in[0,T]$, any $\varepsilon>0$
$$
|u^{\varepsilon}(t)-u^{\varepsilon}(\tau)|^2\lesssim |G(u^{\varepsilon}(t))-G(u^{\varepsilon}(\tau))|,
$$
thus we have that for $\tilde{P}-a.s. \omega$
$$
\|u^{\varepsilon}(t)-u^{\varepsilon}(\tau)\|_{L^2}^2\lesssim \|G(u^\varepsilon(t))-G(u^\varepsilon(\tau))\|_{L^1}.
$$
Since $G(u^\varepsilon)\to g$ in $C^\beta([0,T],L^1)$ for any $\beta<\frac{1}{12}$, let $\varepsilon\to 0$ we have that
$$
\limsup_{\varepsilon\to 0}\|u^{\varepsilon}(t)-u^{\varepsilon}(\tau)\|_{L^2}^2\lesssim \|g(t)-g(\tau)\|_{L^1}\lesssim|t-\tau|^{\beta}.
$$
Since $u^\varepsilon(t)\to u(t)$ in $L^2_w$, by the weakly lower-semicontiniuty, we have that
$$
\|u(t)-u(s)\|_{L^2}^2\leq\liminf_{\varepsilon\to 0}\|u^{\varepsilon}(t)-u^{\varepsilon}(\tau)\|_{L^2}\lesssim |t-\tau|^\beta.
$$
Hence we obtain that $u\in C^{\frac{\beta}{2}}([0,T];L^2)$ $\mathbb{P}-a,s,$. Similarly we have for any $\varepsilon,h>0$
$$
\|u^{\varepsilon}-u^h\|_{L^2}^2\lesssim \|G(u^\varepsilon)-G(u^h)\|_{L^1},\quad \tilde{\mathbb{P}}-a.s..
$$
Let $h\to 0$, we obtain
$$
\|u^{\varepsilon}-u\|_{L^2}^2\lesssim \|G(u^\varepsilon)-g\|_{L^1},\quad \tilde{\mathbb{P}}-a.s.,
$$
which implies that $u^\varepsilon\to u$ in $C^{\frac{\beta}{2}}([0,T];L^2)$ $\tilde{\mathbb{P}}-a.s.$. 

On the other hand, by (\ref{3.6}) we know that 
$$
\mathbb{E}\sup_{t\in[0,T]}\int_\mathcal{D}(|u^\varepsilon|-1)^2dx\lesssim\mathbb{E}\sup_{t\in[0,T]}\|F(u^\varepsilon)\|_{L^1}\lesssim\varepsilon.
$$
As $\varepsilon\to 0$, we have that for $\tilde{\mathbb{P}}-a.s.$ $|u|\equiv 1$ in $L^2$, $\forall t\in[0,T]$, which implies that  for $\tilde{\mathbb{P}}-a.s.$ there exists a measurable set $E(\omega)$ in $\mathcal{D}_T$, such that
$$
u=-1+2\chi_E,\quad\tilde{\mathbb{P}}-a.s..
$$
Since $u$ is $\{\tilde{\mathcal{F}}_t\}_{t\in[0,T]}$-adapted in $L^2$, we know $\chi_E$ is also $\{\tilde{\mathcal{F}}_t\}_{t\in[0,T]}$-adapted in $L^2$.

Moreover by the right hand side of (\ref{3.10}), we obtain that for $\tilde{\mathbb{P}}-a.s.\omega$,
$$
\|g-G(u)\|_{L^1}=\lim_{\varepsilon\to 0}\|G(u^\varepsilon)-G(u)\|_{L^1}\lesssim\lim_{\varepsilon\to 0}\|u^\varepsilon-u\|=0
$$
which implies that $g=G(u)=2S\chi_E$.
Hence we proved the assertion (\rmnum{3}).

Using the estimate (\ref{3.14}), we have that for any $t,\tau\in[0,T]$,
$$
\mathbb{E}\|\chi_{E_t}-\chi_{E_\tau}\|_{L^1}^{2p}\lesssim\mathbb{E}\|\chi_{E_t}-\chi_{E_\tau}\|_{L^2}^{2p}\lesssim\lim_{\varepsilon\to 0}\mathbb{E}\|u^\varepsilon(t)-u^\varepsilon(\tau)\|_{L^2}^{2p}\lesssim|t-\tau|^{\frac{p}{6}}.
$$
Then the assertion (\rmnum{4}) followed by the Kolmogorov's criteria.

Note that the equation (\ref{1.1}) is conserved, i.e. for any $t\in[0,T]$,
$$
\int_\mathcal{D}u^\varepsilon(t,x)dx\equiv \int_\mathcal{D}u^\varepsilon_0(x)dx=|\mathcal{D}|m_0.
$$
Since $u^\varepsilon\to u=-1+2\chi_E$, we have that $|E_t|=\frac{1+m_0}{2}|\mathcal{D}|$. This proved the assertion (\rmnum{5}).

Finally set $g^\varepsilon:=G(u^\varepsilon)$, by (\ref{3.9}) we know that 
$$
|Dg^\varepsilon(t,\cdot)|(\mathcal{D})=\int_{\mathcal{D}}|\nabla g^\varepsilon(t,x)|dx \leq\mathcal{E}^\varepsilon(t).
$$
As $\varepsilon\searrow 0$, since $g^\varepsilon\to g=2S\chi_E$ in $C([0,T];L^1)$ and $\mathcal{E}^\varepsilon\to\mathcal{E}$ in $L^\infty_{w*}(0,T)$, by \cite[Proposition 3.13]{Ambrosio:2000wf}, we obtain that $Dg^\varepsilon\to Dg$ in $L^\infty(0,T;BV)$. Then
by the lower semicontiniuty of the $BV$ norm we obtain that $|D\chi_{E_t}|(\mathcal{D})=\frac{1}{2S}|Dg(t,\cdot)|\leq\frac{1}{2S}\mathcal{E}(t)$. This completes the proof of the theorem.

$\hfill\Box$
\vskip.10in

\subsection{Proof of Theorem \ref{t2.3}}\label{ss3.5}

Now we are in a position to prove Theorem \ref{t2.3}. Before we begin the proof, we need to first recall some crucial lemmas to estimate the following "discrepancy" measure $\zeta^\varepsilon(u^\varepsilon)dx$
\begin{equation}\label{dm}
\zeta^\varepsilon(u^\varepsilon)dx:=\left(\frac{\varepsilon}{2}|\nabla u^\varepsilon|^2-\frac{1}{\varepsilon}F(u^\varepsilon)\right)dx.
\end{equation}

\begin{lemma}(\cite[Lemma 4.4, Theorem 3.6]{Chen:1996bg}).\label{ldm}
Let 
$$
\mathcal{K}^\varepsilon:=\left\{(u,v)\in H^2(\mathcal{D})\times L^2(\mathcal{D}):v=-\varepsilon\Delta u+\frac{1}{\varepsilon}f(u)\;\;\text{in}\;\;\mathcal{D},\;\;\frac{\partial u}{\partial n}=0\;\;\text{on}\;\;\partial \mathcal{D}\right\}.
$$
There exist psitive constants $C_0$ and $\eta_0\in(0,1]$
such that for every $\eta\in(0,\eta_0]$, every $\varepsilon\in(0,1]$, and every $(u^\varepsilon,v^\varepsilon)\in\mathcal{K}^\varepsilon$,
\begin{equation}\label{dm1}
\int_{\left\{x \in \mathcal{D} ;\left|u^{\varepsilon}\right| \geq 1-\eta\right\}}\left[e^{\varepsilon}\left(u^{\varepsilon}\right)+\varepsilon^{-1} f^{2}\left(u^{\varepsilon}\right)\right] \leq C_{0} \eta \int_{\left\{x \in \mathcal{D} ;\left|u^{\varepsilon}\right| \leq 1-\eta\right\}} \varepsilon\left|\nabla u^{\varepsilon}\right|^{2}+C_{0} \varepsilon \int_{\mathcal{D}}( v^{\varepsilon})^2.
\end{equation}
Moreover there exist continuous, non-increasing, and positive functions $M_1(\eta)$ and $M_2(\eta)$ defined on $(0,\eta_0]$ such that for every $\eta\in(0,\eta_0]$, every $\varepsilon\in(0,\frac{1}{M_1(\eta_0)}]$, and every $(u^\varepsilon,v^\varepsilon)\in\mathcal{K}^\varepsilon$, we have that
\begin{equation}\label{dm2}
\int_\mathcal{D}\left(\zeta^\varepsilon(u^\varepsilon)\right)^+dx\leq\eta\int_\mathcal{D}e^\varepsilon(u^\varepsilon)dx+\varepsilon M_2(\eta)\int_\mathcal{D}(v^\varepsilon(x))^2dx,
\end{equation}
where $\left(\zeta^\varepsilon(u^\varepsilon)\right)^+$ is the positive part of $\zeta^\varepsilon(u^\varepsilon)$.
\end{lemma}

\vskip.10in
\textbf{Proof of Theorem \ref{t2.3}}
\vskip.10in
Let $\{u_0^\varepsilon(\cdot)\}_\varepsilon$ be a family of initial data satisfying (\ref{1.2}).  Let $(u^\varepsilon,v^\varepsilon)$ be the solution of (\ref{1.1}) with initial value $u_0^\varepsilon$. The first three assertions can be obtained directly by Theorem \ref{t3.5}. 

In the following we fixed $\omega$ such that all the assertions in Theorem \ref{t3.5} hold. For simplicity of notation, we also denote $\varepsilon_k$ by $\varepsilon$ and omit the notation tilde $\;\tilde{}\;$ in the Theorem \ref{t3.5}.

Since $G(u^\varepsilon)\to 2S\chi_E$ and $|DG(u^\varepsilon)|\leq e^\varepsilon(u^\varepsilon)$ for every $\varepsilon$ and every $(t,x)\in\mathcal{D}_T$, by the lower semicontinuity of the BV norms, we have that
\begin{equation*}
2S|D\chi_{E_t}|dtdx\leq d\mu,
\end{equation*}
which is the inequality (\ref{2.6}).

For any $\psi\in C_c^1([0,t)\times\bar{\mathcal{D}})$, denote $h(t,u):=\int_\mathcal{D}(1+u(x))\psi(t,x)dx$. Since $(u^\varepsilon,v^\varepsilon)$ is a solution to equation (\ref{1.1}), by It\^o's formula we have that for any $\tau\in(0,t)$
$$
h(t,u^\varepsilon(t))-h(0,u^\varepsilon(0))=\int_0^t\int_\mathcal{D}\partial_t\psi(\tau,x)(1+u^\varepsilon(\tau,x))dxd\tau+\int_0^t\langle\psi(\tau,\cdot),du^\varepsilon(\tau)\rangle,
$$
combined with $\psi(t)\equiv 0$, which yields that
\begin{align*}
-\int_\mathcal{D}(1+u^\varepsilon(0,x))\psi(0,x)dx=&\int_0^t\int_\mathcal{D}\partial_t\psi(\tau,x)(1+u^\varepsilon(\tau,x))dxd\tau-\int_0^t\int_\mathcal{D}\nabla v^\varepsilon\nabla\psi\\
&+\varepsilon^\sigma\int_0^t\langle\psi(\tau,\cdot),dW_\tau\rangle.
\end{align*}
Let $\varepsilon\searrow 0$, we obtain that the identity (\ref{2.8}).

In addition, for any $t\in(0,T]$, $\vec{Y}\in C_c^1(\mathcal{D}_t;\mathbb{R}^d)$, a direct calculation by integration by parts yields that
\begin{align*}
\int_\mathcal{D}\vec{Y}\cdot\nabla u^\varepsilon v^\varepsilon&=\int_\mathcal{D}\vec{Y}\cdot\nabla u^\varepsilon\left(-\varepsilon\Delta u^\varepsilon+\frac{1}{\varepsilon}f(u^\varepsilon)\right)\\
&=-\int_\mathcal{D}D\vec{Y}:\left(e^\varepsilon(u^\varepsilon)-\varepsilon\nabla u^\varepsilon\otimes\nabla u^\varepsilon\right)+\int_{\partial \mathcal{D}}e^\varepsilon(u^\varepsilon)\vec{Y}\cdot\vec{n}_{\partial \mathcal{D}}\\
&=-\int_\mathcal{D}D\vec{Y}:\left(e^\varepsilon(u^\varepsilon)-\varepsilon\nabla u^\varepsilon\otimes\nabla  u^\varepsilon\right).
\end{align*}
The last equality holds because $\mathcal{D}$ is an open domain thus $\vec{Y}\equiv 0$ on $\partial\mathcal{D}$.
Then taking integration from $s=0$ to $s=t$ and letting $\varepsilon\searrow 0$, we obtain
\begin{equation}\label{3.20}
\int_0^t2\chi_E\text{div}(v\vec{Y})dxds=\int_0^t\int_\mathcal{D}D\vec{Y}:\left(\text{I} d\mu-(d\mu_{ij})_{d\times d}\right).
\end{equation}

It remains to construct $V$ to finish the proof.
Note that for any $0<\tau<t<T$,
\begin{equation}\label{3.21}
\int_\tau^t\int_{\bar{\mathcal{D}}}d\mu(s,x)=\lim_{\varepsilon\searrow 0}\int_\tau^t\int_{{\mathcal{D}}}e^\varepsilon(u^\varepsilon)dsdx=\int_\tau^t\mathcal{E}(s)ds.
\end{equation}
Therefore, in the sense of Radon measure,
$$
d\mu(t,x)=d\mu^t(x)dt.
$$
By (\ref{3.21}) we have $\mu^t(\bar{\mathcal{D}})=\mathcal{E}(t)$ for $a.e.$ $t\in (0,T]$. Consequetly, for $a.e.$ $t\in (0,T]$ and $a.e.$ $\tau\in (0,t)$, by (\ref{3.4}), we have that
\begin{align*}
\mu^t(\bar{\mathcal{D}})=&\mathcal{E}(t)=\lim_{\varepsilon\searrow 0}\mathcal{E}^\varepsilon(t)=\lim_{\varepsilon\searrow 0}\left(\mathcal{E}^\varepsilon(u^\varepsilon)(\tau)-\int_\tau^t\int_\mathcal{D}|\nabla v^\varepsilon|^2\right)\\
&+\lim_{\varepsilon\searrow 0}\left(\varepsilon^{2\sigma+1}\int_\tau^t\text{Tr}(-\Delta Q)ds+\frac{\varepsilon^{2\sigma-1}}{2}\int_\tau^t\text{Tr}(f'(u^\varepsilon)Q)ds+\varepsilon^\sigma \int_\tau^t\langle v^\varepsilon,dW_\tau\rangle\right).
\end{align*}
Similar as in the proof of Lemma \ref{l3.1}, we have that for $\sigma>\frac{1}{2}$,
$$
\lim_{\varepsilon\searrow 0}\left(\varepsilon^{2\sigma+1}\int_\tau^t\text{Tr}(-\Delta Q)ds+\frac{\varepsilon^{2\sigma-1}}{2}\int_\tau^t\text{Tr}(f'(u^\varepsilon)Q)ds+\varepsilon^\sigma \int_\tau^t\langle v^\varepsilon,dW_\tau\rangle\right)=0.
$$
Hence we deduce that
$$
\mu^t(\bar{\mathcal{D}})\leq\mathcal{E}(\tau)-\int_\tau^t\int_\mathcal{D}|\nabla v|^2dxds=\mu^\tau(\bar{\mathcal{D}})-\int_\tau^t\int_\mathcal{D}|\nabla v|^2dxds,
$$
which is the inequality (\ref{2.10}).

Next, we study the relation between $\mu_{ij}$ and $\mu$. Observe that for any $t\in(0,T]$, and $\vec{Y},\vec{Z}\in C(\bar{\mathcal{D}}_t;\mathbb{R}^d)$,
\begin{equation}\label{3.22}
\begin{aligned}
\varepsilon\int_0^t\int_\mathcal{D} \vec{Y}^T\left(\nabla u^\varepsilon\otimes\nabla u^\varepsilon\right)\vec{Z}
=&\varepsilon\int_0^t\int_\mathcal{D}\sum_{i,j}Y^iZ^j\partial_{x_i}u^\varepsilon\partial_{x_j}u^\varepsilon dxdt\\
\leq&  \varepsilon\int_0^t\int_\mathcal{D}|\vec{Y}||\vec{Z}||\nabla u^\varepsilon|^2dxdt\\
\leq& \int_0^t\int_\mathcal{D}|\vec{Y}||\vec{Z}|e^\varepsilon(u^\varepsilon)+\int_0^t\int_\mathcal{D}|\vec{Y}||\vec{Z}|\zeta^\varepsilon(u^\varepsilon),
\end{aligned}
\end{equation}
where $\vec{Y}^T$ is the transpose of vector $\vec{Y}$. Here in the last inequality we use the definition of $\zeta^\varepsilon(u^\varepsilon)$ in (\ref{dm}) then $e^\varepsilon(u^\varepsilon)+\zeta^\varepsilon(u^\varepsilon)=\varepsilon|\nabla u^\varepsilon|^2$.

By taking $\eta$  as small as enough in (\ref{dm2}), we have that
$$
\lim_{\varepsilon\searrow 0}\int_0^t\int_\mathcal{D}|\vec{Y}||\vec{Z}|\zeta^\varepsilon(u^\varepsilon)\leq 0.
$$
Thus letting $\varepsilon\searrow 0$ in (\ref{3.22}), we obtain that
\begin{equation}\label{3.23}
\int_0^t\int_{\bar{\mathcal{D}}}\vec{Y}^T\left(d\mu_{ij}\right)_{d\times d}\vec{Z}\leq \int_0^t\int_\mathcal{D}|\vec{Y}||\vec{Z}|d\mu.
\end{equation}
Therefore, in the sense of measure $|d\mu_{ij}(t,x)|\leq d\mu(t,x)$. Consequently, there exists $\mu$-measurable functions $\nu_{ij}(t,x)$ such that
$$
d\mu_{ij}(t,x)=\nu_{ij}(t,x)d\mu(t,x),\quad\mu-a.e.\;\;(t,x)\in\bar{\mathcal{D}}_T.
$$
By the definition of $\mu_{ij}$ and (\ref{3.23}), we have that 
$$
0\leq\left(\nu_{ij}\right)_{d\times d}=\left(\nu_{ij}(t,x)\right)_{d\times d}\leq \text{I},\quad\mu-a.e.\;\;(t,x)\in\bar{\mathcal{D}}_T.
$$
Therefore we have that
$$
\left(\nu_{ij}\right)_{d\times d}=\sum_{i=1}^d\lambda_i\vec{\nu}_i\otimes\vec{\nu}_i,\quad\mu-a.e.,
$$
where $\vec{\nu}_i$,  $i=1,\cdots,d$ are $\mu$-measurable unit vectors and $\lambda_i$, $i=1\cdots,d$ are $\mu$-measurable functions, which satisfy
\begin{equation}
0\leq\lambda_i\leq 1\;(i=1,\cdots,d),\quad\sum_{i=1}^d\lambda_i\leq 1,\quad\sum_{i=1}^d\vec{\nu}_i\otimes\vec{\nu}_i=\text{I},\quad\mu-a.e..
\end{equation}
It then follows from (\ref{3.20}) that for $a.e.\;t\in(0,T]$ and for every $\vec{Y}\in C_c^1(\mathcal{D},\mathbb{R}^d)$,
\begin{align*}
2\int_\mathcal{D}\chi_{E_t}\text{div}\left(v(t,x)\vec{Y}(x)\right) dx
&=\int_\mathcal{D}\nabla\vec{Y}(x):\left(\text{I}-\sum_{i=1}^d\lambda_i(t,x)\vec{\nu}_i(t,x)\otimes\vec{\nu}_i(t,x)\right)d\mu^t(x)\\
&=\int_\mathcal{D}\nabla\vec{Y}(x):\sum_{i=1}^dc_i^t(x)\left(\text{I}-\vec{\nu}_i(t,x)\otimes\vec{\nu}_i(t,x)\right)d\mu^t(x),
\end{align*}
where
$$
c_i^t(x)=\lambda_i(t,x)+\frac{1}{d-1}\left(1-\sum_{i=1}^d\lambda_i(t,x)\right).
$$
Clearly, for $a.e.\;t\in(0,T]$, $0\leq c_i^t\leq 1$ and $\sum_{i=1}^dc_i^t\geq 1$ for $\mu^t-a.e.$. Define $p_i^t=\{\vec{\nu}_i(t,x),-\vec{\nu}_i(t,x)\}\in P$ and $V^t$ as in (\ref{2.7}), then $V$ is defined by $dV(t,x,p)=dV^t(x,p)dt$, i.e.
$$
dV(t,x,p)=\sum_{i=1}^dc_i^t(x)\delta_{p_i^t(x)}(p)d\mu^t(x)dpdt,
$$
satisfying (\rmnum{3}) of Definition \ref{d2.5}.

Then by (\ref{2.3}), 
$$
\int_\mathcal{D}\chi_{E_t}(x)\mathrm{div}\left(v(t,x)\vec{Y}(x)\right)dx=\frac{1}{2}\langle\delta V^t,\vec{Y}\rangle.
$$
Thus we obtain (\ref{2.9}). Hence we proved
(\rmnum{4}) of Theorem \ref{t2.3}. This completes the proof of Theorem \ref{t2.3}.

\subsection{The case that $\sigma=\frac{1}{2}$}

As what is shown in the last subsection, for $\sigma=\frac{1}{2}$, the limit of solution to equation (\ref{1.1}) satisfies all the defintion in Definition \ref{d2.1} except (\ref{2.10}). Instead we have
\begin{proposition} \label{p3.8}
Let $\mu^t$ be as in Theorem \ref{t2.3}, then
\begin{equation}\label{5.1}
\mu^t(\bar{\mathcal{D}})+\int_\tau^t\int_\mathcal{D}|\nabla v|^2\leq\mu^\tau(\bar{\mathcal{D}})+C_Q(t-\tau),\;\;\tilde{P}-a.s.,
\end{equation}
where $C_Q:=\text{Tr}(Q)$. 
\end{proposition}
\proof
By using the method as in subsection \ref{ss3.5}, we have that for $\sigma=\frac{1}{2}$
\begin{equation*}
\begin{aligned}
\mu^t(\bar{\mathcal{D}})=&\mathcal{E}(t)
=\lim_{\varepsilon\searrow 0}\left(\mathcal{E}^\varepsilon(u^\varepsilon)(\tau)-\int_\tau^t\int_\mathcal{D}|\nabla v^\varepsilon|^2\right)\\
&+\lim_{\varepsilon\searrow 0}\left(\varepsilon^2\int_\tau^t\text{Tr}(-\Delta Q)ds+\frac{1}{2}\int_\tau^t\text{Tr}(f'(u^\varepsilon)Q)ds+\varepsilon^{\frac{1}{2}}\langle \int_\tau^tv^\varepsilon,dW_\tau\rangle\right)\\
&\leq\mathcal{E}(\tau)-\int_\tau^t\int_\mathcal{D}|\nabla v|^2dxds+(t-\tau)\text{Tr}(Q).
\end{aligned}
\end{equation*}
The last inequality holds because $|u^\varepsilon|\to 1$ and $f'(u^\varepsilon)=3(u^\varepsilon)^2-1$. Thus we obtain (\ref{5.1}).

$\hfill\Box$
\vskip.10in
\begin{remark}
By Propostion \ref{p3.8} and the analysis in the proof of Theorem \ref{t2.3} in Subsection \ref{ss3.5}.
In the case that $\sigma=\frac{1}{2}$, the energy $\mu^t$ may grow a little faster than that in deterministic case. But as what we will show in the next section, at least in radial symmetric case, the pertubation by the noise $\varepsilon^{\frac{1}{2}}dW$ is not strong enough, such that the limit of equation (\ref{1.1}) also converges to deterministic Hele-Shaw model (in a weak sense). Thus we conjucture that in general for $\mathbb{P}-a.s.\;\omega$, the sharp interface limit of (\ref{1.1}) satisfies the deterministic Hele-Shaw model (\ref{1.3}):
\begin{equation*}
\left\{
   \begin{aligned}
   \Delta v&=0 \;\text{in}\;\mathcal{D}\setminus\Gamma_t,\; t>0,\\
\frac{\partial v}{\partial n}&=0\;\text{on}\;\partial\mathcal{D},\\
   v&=SH\;\text{on}\;\Gamma_t,\\
  \mathcal{V}&=\frac{1}{2}(\partial_{n}v^+-\partial_{n}v^-)\;\text{on}\;\Gamma_t.\\
   \end{aligned}
   \right.
\end{equation*}
\end{remark}

\section{Case of radial symmetry for $\sigma\geq\frac{1}{2}$}\label{s4}

In this section we are going to prove Theorem \ref{t2.8}. In this case of radial symmetry, we assume $\mathcal{D}=B_1$. We denote by $B_r$ the ball centered at the origin $0$ with radius $r$.
We also denote by $S_r$ the sphere of radius $r$ in $\mathbb{R}^d$ and by $\omega_d$ the area of unit sphere $S_1$. 
Any function $u$ in this section of the form $u(x)\equiv u(|x|)$. For convenience, we do not distinguish functions of $x\in B_1$ from functions of $r\in[0,1)$. We only distinguish the integrals of $dx$ from that of $dr$, due to consideration of singularities at the origin.

 Denote $r=|x|$, then the equation (\ref{1.1}) should be changed as
\begin{equation}\label{4.1}
\left\{
\begin{aligned}
&du^{\varepsilon}=\partial_{rr} v^\varepsilon dt+\frac{d-1}{r}\partial_r v^\varepsilon dt+\varepsilon^\sigma dW_t, \quad (t,r)\in[0,T]\times[0,1],\\
&v^\varepsilon=-\varepsilon\partial_{rr} u^{\varepsilon}(t)-\frac{d-1}{r}\partial_r u^\varepsilon+\frac{1}{\varepsilon} f(u^{\varepsilon}(t)),\quad (t,x)\in[0,T]\times[0,1],\\
&\partial_r u^\varepsilon(t,1)=\partial_r v^\varepsilon(t,1)=0, \quad t\in[0,T],\\
&u^{\varepsilon}(0,r)=u_0^\varepsilon(r),\quad x\in\mathcal{D}.
\end{aligned}\right.
\end{equation}
Here $W_t$ is given by 
\begin{equation}\label{4.2}
W_t=\sum_{k\in\mathbb{Z}^d}\alpha_k b_k(r)\beta_k(t),
\end{equation}
where
$\{\beta_k\}_{k\in\mathbb{Z}^d}$ is a sequence of independent Brownian motions and $\{\alpha_k\}$ satisfies
\begin{equation}\label{4.3}
\sum_{k\in\mathbb{Z}^d}\lambda_k\alpha_k^2<\infty.
\end{equation}
$\{b_k\}_{k\in\mathbb{Z}^d}$ is an orthogonal basis in $L_0^2(0,1)$, which is defined as $\{f\in L^2(0,1):\int_0^1 f(r)r^{d-1}dr=0\}$, i.e.
\begin{equation}\label{4.4}
\int_0^1b_k(r)r^{d-1}dr=0,\quad\forall k\in\mathbb{Z}^d.
\end{equation}
Note that (\ref{4.2})-(\ref{4.4}) is just the radially symmetric version of condition (\ref{2.1}) and (\ref{2.2}). Moreover, all the results we obtained in the previous section also hold for this case.
In particular, there exists a Borel set $\tilde{E}\subset[0,T]\times[0,1)$ such that $E=\{(t,x)\in\mathcal{D}_T:(t,|x|)\in\tilde{E}\}$ and $\tilde{E}_t:=\{r\in[0,1):(t,r)\in[0,T]\times[0,1)\}$ is a BV set in $[0,1)$ for any $t\in[0,T]$. 

\begin{remark}
For the existence of radial symmetric solution to (\ref{1.1}) under the assumption in this section, we only need to check that any solution $u^\varepsilon$ to (\ref{1.1}) is invariant under the rotation transformation. Then by the uniqueness, we can obtain that $u^\varepsilon$ is radial symmetric.

In fact, any rotation transformation in $\mathbb{R}^d$ can be indentified as an orthogonal matrix with determinant $1$, i.e. an element in $SO(d)$. For any $A\in SO(d)$, a direct calculation yields that
$$
\nabla (v\circ A)=(\nabla v)\circ A\;\;(\Delta v)\circ A=\Delta(v\circ A).
$$
Then we have that for any solution $(u^\varepsilon,v^\varepsilon)$ to (\ref{1.1}), $(u^\varepsilon\circ A,v^\varepsilon\circ A)$ is also a solution to (\ref{1.1}). By the  uniqueness of solutions to (\ref{1.1}), if the initial value of  $(u^\varepsilon,v^\varepsilon)$  is radial symmetric, $(u^\varepsilon,v^\varepsilon)$ is also radial symmetric. In this case, equation (\ref{1.1}) is equivalent to (\ref{4.1}).

\end{remark}

We also mention that all the results in \cite[Section 5]{Chen:1996bg} only depend on the second equation in (\ref{4.1}) and the estimate of $(u^\varepsilon,v^\varepsilon)$. Thus with a similar proof, we obtain the following theorems.

\begin{theorem}
Assume that $\{(\tilde{u}^\varepsilon,\tilde{v}^\varepsilon)\}$ is obtained in Theorem \ref{t3.5}. Then 
$$
\lim_{\varepsilon\searrow 0}\int_0^T\int_\mathcal{D}\big|\zeta^\varepsilon(\tilde{u}^\varepsilon)\big|dxdt=0,\quad\mathbb{\tilde{P}}-a.s.,
$$
where $\zeta^\varepsilon(\tilde{u}^\varepsilon)$ is the discrepancy measure defined in (\ref{dm}).
\end{theorem}

\proof In this proof, we ignore the notation tilde $\;\tilde{}\;$ for simplicity.

For a fixed $\omega$ such that all the assertions in Theorem \ref{t3.5} hold. By the same proof as \cite[Theorem 5.1]{Chen:1996bg}, we have that there exists a cnostant $C>0$ he following estimates
\begin{equation}\label{4.5}
\int_{B_\delta}e^\varepsilon(u^\varepsilon)dx\leq C\delta M^\varepsilon(t),\;\;\forall\delta\in(0,1),
\end{equation}
where $M^\varepsilon(t):=1+\mathcal{E}^\varepsilon({u}^\varepsilon)(t)+\|{v}^\varepsilon\|_{H^1}^2\in L^1(0,T)$, and 
\begin{equation}\label{4.6}
\sup_{0<r<1}\left|r^{d-1}\left(\zeta^\varepsilon(u^\varepsilon)+v^\varepsilon u^\varepsilon\right)\right|\leq CM^\varepsilon(t).
\end{equation}
Hence for any small $\delta$ and $\eta$,
\begin{align*}
\int_{\mathcal{D}}\left|\zeta^{\varepsilon}\left(u^{\varepsilon}\right)\right| dx
&\leq \int_{B_{\delta} \cup\left\{\left|u^{\varepsilon}\right| \geq 1-\eta\right\}} \left|\zeta^{\varepsilon}\left(u^{\varepsilon}\right)\right|  d x+\int_{\mathcal{D} \cap\left\{r>\delta,\left|u^{\epsilon}\right| \leq 1-\eta\right\}}\left|\zeta^{\varepsilon}\left(u^{\varepsilon}\right)\right| dx\\
&\leq  \int_{B_{\delta} \cup\left\{\left|u^{\varepsilon}\right| \geq 1-\eta\right\}} e^{\varepsilon}\left(u^{\varepsilon}\right) d x+\int_{\mathcal{D} \cap\left\{r>\delta,\left|u^{\epsilon}\right| \leq 1-\eta\right\}}\left[\left|v^{\varepsilon}\right||u^\varepsilon|+r^{1-d} C M^{\varepsilon}(t)\right]dx\\
&\leq \int_{B_{\delta} \cup\left\{\left|u^{\varepsilon}\right| \geq 1-\eta\right\}} e^{\varepsilon}\left(u^{\varepsilon}\right) d x+\int_{\mathcal{D} \cap\left\{r>\delta,\left|u^{\epsilon}\right| \leq 1-\eta\right\}}\left[\left|v^{\varepsilon}\right|(1-\eta)+r^{1-d} C M^{\varepsilon}(t)\right]dx,
\end{align*}
where we used the definition of $\zeta^\varepsilon$ and $e^\varepsilon(u^\varepsilon)$ and (\ref{4.5}) in the sencond inequality.

For the first intergral above, we have that
\begin{align*}
\int_{B_{\delta} \cup\left\{\left|u^{\varepsilon}\right| \geq 1-\eta\right\}} e^{\varepsilon}\left(u^{\varepsilon}\right) d x
&\leq \int_{B_{\delta}} e^{\varepsilon}\left(u^{\varepsilon}\right) d x+\int_{\left\{\left|u^{\varepsilon}\right| \geq 1-\eta\right\}} e^{\varepsilon}\left(u^{\varepsilon}\right) d x\\
&\leq C\delta M^\varepsilon(t)+C_0\eta M^\varepsilon(t)+C_0\varepsilon M^\varepsilon(t),
\end{align*}
where we used (\ref{4.5}) and (\ref{dm1}) in the second inequality.

For the second integral, we have that
\begin{align*}
\int_{\mathcal{D} \cap\left\{r>\delta,\left|u^{\epsilon}\right| \leq 1-\eta\right\}}\left[\left|v^{\varepsilon}\right|(1-\eta)+r^{1-d} C M^{\varepsilon}(t)\right]dx
&\leq \int_{\left\{\left|u^{\epsilon}\right| \leq 1-\eta\right\}}\left[\left|v^{\varepsilon}\right|+\delta^{1-d} C M^{\varepsilon}(t)\right]dx\\
&\leq \mathcal{H}^d(\left\{\left|u^{\epsilon}\right| \leq 1-\eta\right\})\left(M^\varepsilon(t)^{\frac{1}{2}}+\delta^{1-d}CM^\varepsilon(t)\right).
\end{align*}
By Theorem \ref{t3.5} we know that
$
\varepsilon^{-1}F(u^\varepsilon)$ is uniformly bounded in $L^\infty(0,T;L^1)$.
Thus there exists a constant $C_1>0$ such that
$$
\mathcal{H}^d(\left\{\left|u^{\epsilon}\right| \leq 1-\eta\right\})\leq\mathcal{H}^d(\left\{\left|\left|u^{\epsilon}\right|-1\right| \geq \eta\right\})\leq\eta^{-2}\int_\mathcal{D}F(u^\varepsilon)dx\leq C_1\eta^{-2}\varepsilon.
$$

Combining all the estimates above,
we have that for any $\eta,\delta>0$, there exists a constant $C(\delta,\eta)>0$, such that
$$
\int_\mathcal{D}\big|\zeta^\varepsilon(\tilde{u}^\varepsilon)\big|dx\leq C_2\left(\delta+\eta+\varepsilon+C(\delta,\eta)\varepsilon\right)M^\varepsilon(t),
$$
 $C_2$ is independent of $\varepsilon,\eta,\delta$. Integrating the last inequality in $(0,T)$ and letting first $\varepsilon\to 0$ and then $\delta,\eta$ to $0$, we can obtain the theorem.

$\hfill\Box$
\vskip.10in

In the following, we are going to prove
$$
d\mu=2S|D\chi_{E}|dxdt.
$$

\begin{theorem}\label{t4.2}
Let $\{(\tilde{u}^{\varepsilon_k},\tilde{v}^{\varepsilon_k})\}_k$ are radially symmetric solutions of (\ref{4.1}) which satisfy all the assertions in Theorem \ref{t3.5}. Then for any $t\in(0,T]$, $\psi\in C_c(\mathcal{D}_t)$,
$$
\int_0^t\int_\mathcal{D}\psi(t,x)d\mu(t,x)=2S\int_0^t\int_\mathcal{D}\psi|D\chi_{E_t}|dxdt,\quad\tilde{\mathbb{P}}-a.s..
$$
\end{theorem}
\proof 
The proof is a modification of the proof of \cite[Theorem 5.3]{Chen:1996bg}. The only difference is that in stochasic case,
by (\ref{3.4}), Theorem \ref{t2.8} and Proposition \ref{p3.8}, we know that for all $\sigma\geq\frac{1}{2}$, there exists a $h^\varepsilon\in L^2(\Omega,\mathcal{F},\mathbb{P};L^2(0,T))$ such that
$$
\tilde{\mu}^t_\varepsilon(\mathcal{D}):=\int_\mathcal{D}e^\varepsilon(\tilde{u}^\varepsilon(t,x))dx=\tilde{\mathcal{E}}^\varepsilon(t)\leq \mathcal{E}_0+h^\varepsilon(t)\;\;\mathbb{P}-a.s.,
$$
where for $\mathbb{P}-a.s.\;\omega$, $\{h^\varepsilon(\omega,\cdot)\}_\varepsilon$ is bounded in $L^2(0,T)$, while in deterministic case as in  \cite{Chen:1996bg}, $h^\varepsilon$ is just $0$. Then the rest proof just follows  the proof of \cite[Theorem 5.3]{Chen:1996bg} for a fixed $\omega$. We put it into the Appendix.

$\hfill\Box$
\vskip.10in

\textbf{Proof of Theorem \ref{t2.8}}

The definition of $V$ in (\ref{2.7}) can be written as
$$
dV^t(x,p)=\sum_{i=1}^dc_i^t(x)\delta_{p_i^t(x)}d\mu^t(x)dp.
$$
From (\ref{2.5}) we know that 
$$
d\|V^t\|(x)=\sum_{i=1}^dc_i^t(x)d\mu^t(x)\geq d\mu^t(x).
$$

As what we mentioned in Remark \ref{r2.9}, by \cite[Theorem 14.3]{Simon:1983ti},
for  $a.e.\;t\in[0,T]$
$$
\mu^t=2S|D\chi_{E_t}|=2S\mathcal{H}^{d-1}\lfloor \partial^*E_t.
$$
Here $\mathcal{H}^{d-1}\lfloor \partial^*E_t$ is the $(d-1)$-dimensional Haurdorff measure on $\partial^*E_t$. $\partial^*E_t$ is the so-called reduced boundary of $E_t$ (see  \cite[Section 14]{Simon:1983ti} for details). In our case
$$
\partial^*E_t=\{x\in\mathcal{D}:|\vec{\nu}_{E_t}(x)|=1\}=\mathrm{supp}(|D\chi_{E_t}|).
$$
Moreover, $\partial^*E_t$ is a countably $(d-1)$-rectifiable set (see \cite[Chapter 3]{Simon:1983ti} for details) and
\begin{equation}\label{4.13.0}
\lim_{\rho\searrow 0}\frac{|D\chi_{E_t}|(B_{\rho}(x))}{\rho^{d-1}}=1,\;\;\mathcal{H}^{d-1}-a.e.\;x\in\partial^*E_t,
\end{equation}
where $B_\rho(x)$ is the ball in $\mathbb{R}^d$ with radius $\rho$ and centered at $x$.
Since $\|V^t\|\geq\mu^t=2S|\chi_{E_{t}}|$, by \cite[Theorem 42.4]{Simon:1983ti}, $V^t$ is rectifiable. Thus by the definition of rectifiability and the expression of $V^t$, we have that
$$
dV^t(x,p)=2S|D\chi_{E_t}|dx\delta_{\vec{\nu}_{E_t}(t,x)}(dp)\;\;\text{as Radon measure on}\;\;\mathcal{D}\times P,
$$
i.e.
$$
dV(t,x.p)=dV^t(x.p)dt=d\mu^t(t,x)\delta_{\vec{\nu}_{E_t}(t,x)}(dp).
$$
Hence we conclude that $c_1^t=1$, $c_2^t=\cdots=c_d^t=0$ and $p_1^t=\vec{\nu_{E_t}}$.
Then by the construction of $V$ in subsection \ref{ss3.5}, we have that
$\lambda_1=1$, $\lambda_2=\cdots=\lambda_d=0$ and
$$
\left(d\mu_{ij}\right)_{d\times d}=\vec{\nu}_{E_t}\otimes\vec{\nu}_{E_t}d\mu.
$$

Note that by (\ref{E}), $\operatorname{supp}(D\chi_{\tilde{E}_t})$ is a countable set, thus
$$
d|D\chi_{E_t}(x)|=\omega_d\sum_{r\in \operatorname{supp}(D\chi_{\tilde{E}_t})}\delta_r(|x|)dx,
$$
where $\delta_r$ is the Dirac measure on $\mathbb{R}$.
Then the empty set $\emptyset$ is the only zero $\mathcal{H}^{d-1}$ measurble subset of $\partial^*E_t$. By
 (\ref{4.13.0}),
$$
\sup _{x \in \mathcal{D},1>r>0} \frac{\left|D \chi_{E_{t}}\right|(B_\rho(x))}{\rho^{N-1}}<\infty.
$$
Then by \cite[Theorem 7.1]{Soner:1995gb}, $\langle D\chi_{E_t},v\vec{Y}\rangle$ is a bounded linear funtional w.r.t $\vec{Y}$ over $C_c(\mathcal{D};\mathbb{R}^d)$. Namely, $v$ is $|D\chi_{E_{t}}|$-measurable hence is also $\|V^t\|$-measurable.
Then by (\ref{2.9}),
for any $\vec{Y}\in C_c^1(\mathcal{D},\mathbb{R}^d)$, we have that
$$
-\langle\delta V^t,\vec{Y}\rangle=2\langle D\chi_{E_t},v\vec{Y}\rangle=2\langle\vec{\nu}_{E_t}|D\chi_{E_t}|,v\vec{Y}\rangle=\frac{1}{S}\langle\|V^t\|,v\vec{Y}\cdot\vec{\nu}_{E_t}\rangle.
$$
Hence by the Definiton \ref{d2.5}, we obtain that
$$
S\vec{H}_{V^t}=v\vec{\nu}_{E_t}\;\;\|V^t\|-a.e.,
$$
where $\vec{H}_{V^t}$ is the mean curvature vector of $V^t$ in 
Definiton \ref{d2.5}. This also implies that for $\|V^t\|-a.e.\;x\in\mathcal{D}\setminus\mathrm{supp}(|D\chi_{E_t}|)$, $\vec{H}_{V^t}=0$.
Thus we have that
$$
v=S\vec{H}_{V^t}\cdot\vec{\nu}_{E_t}\;\;\text{on}\;\;\mathrm{supp}(|D\chi_{E_t}|).
$$

\section{The case for "smeared" noise}\label{s6}
We observe that the requirement $\sigma\geq \frac{1}{2}$ only comes from the second variation term $\frac{\varepsilon^{2\sigma-1}}{2}\text{Tr}(f'(u^\varepsilon)Q)$ in (\ref{3.4}) when we apply It\^o's formula on $\mathcal{E}^\varepsilon(u^\varepsilon)$. If there were no such term $\frac{\varepsilon^{2\sigma-1}}{2}\text{Tr}(f'(u^\varepsilon)Q)$, Theorem \ref{t3.5} would hold for $\sigma\geq 0$.

This motivates us to consider the following equation:
\begin{equation}\label{6.1}
\left\{
\begin{aligned}
&\frac{\partial u^{\varepsilon}}{\partial t}=\Delta v^\varepsilon+\varepsilon^\sigma \xi^\varepsilon_t, \quad (t,x)\in[0,T]\times\mathcal{D},\\
&v^\varepsilon=-\varepsilon\Delta u^{\varepsilon}(t)+\frac{1}{\varepsilon} f(u^{\varepsilon}(t)),\quad (t,x)\in[0,T]\times\mathcal{D},\\
&\frac{\partial u^{\varepsilon}}{\partial n}=\frac{\partial v^{\varepsilon}}{\partial n}=0, \quad (t,x)\in[0,T]\times\partial\mathcal{D},\\
&u^{\varepsilon}(0,x)=u_0^\varepsilon(x),\quad x\in\mathcal{D},
\end{aligned}\right.
\end{equation}
where $u_0^\varepsilon$ satisfies (\ref{1.2}) and $\xi^\varepsilon_t$ is formally defined by $\xi^\varepsilon_t=\int_{-\infty}^\infty\rho_\varepsilon(t-s)dW_s$. In fact, let $(W_t,t\geq 0)$ be a $Q$-Wiener process on $L_0^2(\mathcal{D})$ defined on a prabobility space $(\Omega,\mathcal{F},\mathbb{P})$, where $Q$ satisfies (\ref{2.1}) and (\ref{2.2}). We extend the definition of $(W_t,t\geq 0)$ to negative time by considering an i.i.d $Q$-Wiener process $(\hat{W}_t,t\geq 0)$ and setting $W_t=\hat{W}_{-t}$ for $t< 0$. Then $(W_t,t\in\mathbb{R})$ is a two-sided $Q$-Wiener process on $L_0^2$. Let $\rho$ be a mollifying kernel i.e. 
$$
\rho\in C^\infty(\mathbb{R}),\quad 0\leq\rho\leq 1,\quad\text{supp}\rho\subset [-1,1],\quad\int_{\mathbb{R}}\rho=1,\quad\rho(t)=\rho(-t).
$$
For $\gamma>0$ we set $\rho_\varepsilon(t)=\varepsilon^{-\gamma}\rho(\frac{t}{\varepsilon^\gamma})$. Then the approximate Wiener process $W^\varepsilon_t$ is defined as
\begin{equation}\label{6.2}
W^\varepsilon_t:=\int_{-\infty}^\infty\rho_\varepsilon(t-s)W_sds,
\end{equation}
Its derivative is defined as
\begin{equation}\label{6.3}
\xi^\varepsilon_t:=\frac{dW^\varepsilon_t}{dt}=\int_{-\infty}^\infty\rho_\varepsilon(t-s)dW_s.
\end{equation}
Since $\rho_\varepsilon$ is supported on $[-\varepsilon^\gamma,\varepsilon^\gamma]$, only the definition on negative time $[-\varepsilon^\gamma,0)$ of $W_t$ is used.
Thus we have that for any $g\in L^2(\mathcal{D})$
\begin{equation}\label{6.4}
\begin{aligned}
\int_0^T\langle g(t),\xi^\varepsilon_t\rangle dt
&=\int_0^T\langle g(t),\int_{-\varepsilon^\gamma}^{t+\varepsilon^\gamma}\rho_\varepsilon(t-s)dW_s\rangle dt\\
&=\int_{-\varepsilon^\gamma}^{T+\varepsilon^\gamma}\langle \int_0^T\rho_\varepsilon(s-t)g(t)dt,dW_s\rangle.
\end{aligned}
\end{equation}

\begin{lemma}\label{l6.1}
There exists a constant which only depends on $T$ such that for any $\varepsilon\in(0,1]$ and any $p\geq 1$, any $\sigma\geq 0$
\begin{equation}
\mathbb{E}\sup_{t\in[0,T]}\mathcal{E}^\varepsilon(t)^p\leq C_T(\varepsilon^{\sigma}+\mathcal{E}_0)^p, 
\end{equation}
and
\begin{equation}
\mathbb{E}\left(\int_0^T\|\nabla v^\varepsilon\|_{L^2}^2dt\right)^p \leq C_T(\varepsilon^{\sigma}+\mathcal{E}_0)^p.
\end{equation}
\end{lemma}
\proof
The proof is a modification of Lemma \ref{l3.1}. 

Note that the noise in equation \ref{6.1} is smooth in time, which enable us to apply Newton-Leibniz formula on $\mathcal{E}^\varepsilon$ to avoid the second variation term in (\ref{3.4}). We have that
\begin{equation}
\frac{d\mathcal{E}^\varepsilon(u^\varepsilon)}{dt}=\langle D\mathcal{E}^\varepsilon(u^\varepsilon),\partial_t u^\varepsilon\rangle=-\langle\nabla v^\varepsilon,\nabla v^\varepsilon\rangle+\varepsilon^\sigma\langle v^\varepsilon,\xi^\varepsilon_t \rangle.
\end{equation}
By (\ref{6.4}) we know that
$$
\int_0^T\langle v^\varepsilon(t),\xi^\varepsilon_t\rangle dt
=\int_{-\varepsilon^\gamma}^{T+\varepsilon^\gamma}\langle \rho_\varepsilon*v^\varepsilon(t),dW_t\rangle,
$$
where we simiply denote 
$$
\rho_\varepsilon*v^\varepsilon(t):=\int_0^T\rho_\varepsilon(t-s)v^\varepsilon(s)ds.
$$
Similarly as the proof in Lemma \ref{l3.1}. by Burkholder-Davis-Gundy type inequality, we have that
\begin{align*}
\mathbb{E}\sup_{t\in[0,T]}|\int_0^T\langle v^\varepsilon(t),\xi^\varepsilon_t\rangle dt|
&\lesssim
\mathbb{E}\sup_{t\in[0,T]}|\int_0^T\langle \rho_\varepsilon*v^\varepsilon(t),dW_t\rangle|\\
&\lesssim \left(\mathbb{E}\sup_{t\in[0,T]}|\int_0^T\langle \rho_\varepsilon*v^\varepsilon(t),dW_t\rangle|^2\right)^{\frac{1}{2}}\\
&\lesssim \left(\int_0^T\mathbb{E}\|\sqrt{Q}(\rho_\varepsilon*v^\varepsilon(t))\|_{L^2}^2dt\right)^{\frac{1}{2}}\\
&\lesssim \left(\int_0^T\left(\rho_\varepsilon*\mathbb{E}\|\nabla v^\varepsilon\|_{L^2}\right)^2dt\right)^{\frac{1}{2}}\\
&\lesssim \left(\int_0^T\mathbb{E}\|\nabla v^\varepsilon\|_{L^2}^2dt\right)^{\frac{1}{2}},
\end{align*}
where we used the Young's inequality in the last inequality. The rest is the same as in the proof of Lemma \ref{l3.1}. We omit it here for simplicity.

$\hfill\Box$
\vskip.10in

With the same notation and proof as in Theorem \ref{t3.5}, we can obtain a tightness result for any $\sigma\geq 0$.

\begin{theorem}\label{t6.2} Assume $
\sigma\geq 0$, $Q$ satisfies (\ref{2.1}) and (\ref{2.2}).
There exist a probability space $(\tilde{\Omega},\tilde{\mathcal{F}},\tilde{\mathbb{P}})$ on $\mathcal{X}^1\times\mathcal{X}^2$, a subsequence $\varepsilon_k$ (we still denote it as $\varepsilon$ for simplicity) and 
$$\left\{\left(\varepsilon^{-1}\sup_{t\in[0,T]}\|F(\tilde{u}^\varepsilon)\|_{L^1},\mathcal{E}^\varepsilon(\tilde{u}^\varepsilon),\tilde{u}^\varepsilon,G(\tilde{u}^\varepsilon),\tilde{v}^\varepsilon,e^\varepsilon(\tilde{u}^\varepsilon)dxdt,\{\varepsilon\partial_{x_i}\tilde{u}^\varepsilon\partial_{x_j}\tilde{u}^\varepsilon dxdt\}_{ij}\right)\right\}\subset \mathcal{X}^1\times\mathcal{X}^2$$ 
and 
$$\left(a,\mathcal{E},u,g,v,\mu,\{\mu_{ij}\}_{ij}\right)\in \mathcal{X}^1\times\mathcal{X}^2,$$ 
such that

(\rmnum{1}) $\tilde{P}\circ \left(\varepsilon^{-1}\sup_{t\in[0,T]}\|F(\tilde{u}^\varepsilon)\|_{L^1},\mathcal{E}^\varepsilon(\tilde{u}^\varepsilon),\tilde{u}^\varepsilon,G(\tilde{u}^\varepsilon),\tilde{v}^\varepsilon,e^\varepsilon(\tilde{u}^\varepsilon)dxdt,\{\varepsilon\partial_{x_i}\tilde{u}^\varepsilon\partial_{x_j}\tilde{u}^\varepsilon dxdt\}_{ij}\right)^{-1}=\hat{\mathbb{P}}^\varepsilon$ on $\mathcal{X}^1\times\mathcal{X}^2$,

(\rmnum{2}) $\left(\varepsilon^{-1}\sup_{t\in[0,T]}\|F(\tilde{u}^\varepsilon)\|_{L^1},\mathcal{E}^\varepsilon(\tilde{u}^\varepsilon),\tilde{u}^\varepsilon,G(\tilde{u}^\varepsilon),\tilde{v}^\varepsilon,e^\varepsilon(\tilde{u}^\varepsilon)dxdt,\{\varepsilon\partial_{x_i}\tilde{u}^\varepsilon\partial_{x_j}\tilde{u}^\varepsilon dxdt\}_{ij}\right)$ converges to $\left(0,\mathcal{E},u,g,v,\mu,\{\mu_{ij}\}_{ij}\right)$ in $\mathcal{X}^1\times\mathcal{X}^2$, $\tilde{\mathbb{P}}-a.s$, as $\varepsilon\searrow 0$.

In particular, for $\tilde{\mathbb{P}}-a.s. \omega$, there exists a borel set $E\in \mathcal{D}_T$, such that as $\varepsilon\searrow 0$

(\rmnum{3}) $u^\varepsilon\to u$ in $u$ in $C^\beta([0,T];L^2)$, $g=G(u)=2S\chi_E$ a.e. in $\mathcal{D}_T$ and in $C^\beta([0,T];L^1)$, $u=-1+2\chi_E$ a.e. in $\mathcal{D}_T$ and in $C^\beta([0,T];L^2)$.

Moreover, denote $E_t:=\left\{x:(t,x)\in E\right\}$. Then 

(\rmnum{4}) For all $\beta<\frac{1}{12}$, $\tilde{\mathbb{P}}\left(\chi_{E}\in C^\beta([0,T];L^1)\right)=1$,

(\rmnum{5}) $\tilde{\mathbb{P}}\left(|E_t|=|E_0|=\frac{1+m_0}{2}|\mathcal{D}|,\forall t\in[0,T]\right)=1$,

(\rmnum{6}) $\tilde{\mathbb{P}}\left(\chi_{E}\in L^\infty(0,T;BV)\right)=1$.
\end{theorem}
\proof
For all $\sigma\geq 0$,  one can check that with Lemma \ref{l6.1} true,
all the estimate in Subsection \ref{ss3.3} and \ref{ss3.4} hold for the solution $(u^\varepsilon,v^\varepsilon)$ to equation (\ref{6.1}). Then the same proof as Theorem \ref{t3.5} follows.

$\hfill\Box$
\vskip.10in

Moreover, for $\sigma\geq 0$, the same argument as in Subsection \ref{ss3.5} yields that
\begin{theorem}\label{t6.3}
Assume that $\sigma\geq 0$ and (\ref{1.2}) hold. Let $(u^\varepsilon,v^\varepsilon)$ be the solution to (\ref{6.1}). Then there exist a probability space $(\tilde{\Omega},\tilde{\mathcal{F}},\tilde{\mathbb{P}})$, $(\tilde{u}^\varepsilon,\tilde{v}^\varepsilon)\in C([0,T];L^2)\times L^2(0,T;H^1)$ with  $\tilde{\mathbb{P}}\circ\left(\tilde{u}^\varepsilon,\tilde{v}^\varepsilon\right)^{-1}=\mathbb{P}\circ\left(u^\varepsilon,v^\varepsilon\right)^{-1}$ on $C([0,T];L^2)\times L^2(0,T;H^1)$. There also exists a subsequence $\varepsilon_k$ such that as $\varepsilon_k\searrow 0$ the following holds:

(\rmnum{1}) There exists a measurable set $E\subset\tilde{\Omega}\times\mathcal{D}_T$, such that for $\tilde{\mathbb{P}}-a.s.\; \omega$
$$
\tilde{u}^{\varepsilon_k}(\omega)\to -1+2\chi_{E(\omega)},\quad a.e. \;\;\text{in}\;\;\mathcal{D}_T\;\;\text{and in}\;\;C^\beta([0,T];L^2)
$$
for any $\beta<\frac{1}{12}$ where $E(\omega):=\{(t,x)\in\mathcal{D}_T:(\omega,t,x)\in E\}$;

(\rmnum{2}) There exists a random variable  $v\in L_w^2(0,T;H^1)$ ($v$ is weakly measurable in $L^2(0,T;H^1)$)  such that for $\tilde{\mathbb{P}}-a.s.\;\omega$
$$
\tilde{v}^{\varepsilon_k}(\omega)\to v(\omega) \quad\text{weakly in}\;\; L^2(0,T;H^1);
$$

(\rmnum{3}) There exist random variables $\mu\in\mathfrak{M}_R$ and $\{\mu_{ij}\}_{i,j=1}^d\in\mathfrak{M}^{d\times d}$ such that for $\tilde{\mathbb{P}}-a.s.\;\omega$
\begin{equation}
\begin{aligned}
e^{\varepsilon_k(\omega)}(\tilde{u}^{\varepsilon_k(\omega)})dxdt
&\to d\mu(\omega,t,x) \quad\text{weakly in}\;\;\mathfrak{M}_R,\\
\varepsilon_k\partial_{x_i}\tilde{u}^{\varepsilon_k}(\omega)\partial_{x_j}\tilde{u}^{\varepsilon_k}(\omega)dxdt&\to d\mu_{ij}(\omega,t,x)\quad\text{weakly in}\;\;\mathfrak{M},\;\;\forall i,j=1,\cdots,d.
\end{aligned}
\end{equation}

(\rmnum{4}) For $\tilde{\mathbb{P}}-a.s.\;\omega$, there exists a Radon measure $V(\omega)$ on $\mathcal{D}_T\times P$,  and  $\mu^t(\omega,x)dt=d\mu(\omega,t,x)$ such that for any $t\in(0,T]$ and $\vec{Y}\in C_0^1(\mathcal{D}_t;\mathbb{R}^d)$
\begin{equation}
\int_0^t\langle\delta V^s,\vec{Y}\rangle ds=\int_0^t\int_\mathcal{D}\nabla\vec{Y}:\left(Id\mu(s,x)-\left(\mu_{ij}(s,x)\right)_{d\times d}\right).
\end{equation}

In particular,  for $\tilde{\mathbb{P}}-a.s.\;\omega$, $\left(E(\omega),v(\omega),V(\omega)\right)$ satisfies all the properties in Definition \ref{d2.1} except (\ref{2.10}).
 If $\sigma>0$, (\ref{2.10}) holds, thus $\left(E(\omega),v(\omega),V(\omega)\right)$ is a weak solution in the sense of Definition \ref{d2.1}.

\end{theorem}
\proof
The proof is almost the same as in Subsection \ref{ss3.5}. The only difference is in the proof of 
that the exsitence of a $Q$-Wiener process on $L^2$ cannot be otained directly such that for any $\varepsilon>0$, (\ref{3.17.1}) holds. We use the original equation (\ref{6.1}) to prove (\ref{2.8}) directly.

In fact, for any $\psi\in C_c^1([0,t)\times\mathcal{D})$,
\begin{align*}
-\int_\mathcal{D}(1+u^\varepsilon(0,x))\psi(0,x)dx=&\int_0^t\int_\mathcal{D}\partial_t\psi(\tau,x)(1+u^\varepsilon(\tau,x))dxd\tau-\int_0^t\int_\mathcal{D}\nabla v^\varepsilon\nabla\psi dxd\tau\\
&+\varepsilon^\sigma\int_0^t\int_\mathcal{D}\psi(\tau,x)\xi^\varepsilon(\tau,x)dxd\tau.
\end{align*}
Thus for $\mathbb{P}-a.s.\omega\in\Omega$, 
$$
\lim_{\varepsilon\to 0}\left(\int_0^t\int_\mathcal{D}\partial_t\psi(\tau,x)(1+u^\varepsilon(\tau,x))dxd\tau+\int_\mathcal{D}(1+u^\varepsilon(0,x))\psi(0,x)dx-\int_0^t\int_\mathcal{D}\nabla v^\varepsilon\nabla\psi dxd\tau\right)=0.
$$
Since $\tilde{\mathbb{P}}\circ(\tilde{u}^\varepsilon,\tilde{v}^\varepsilon)^{-1}=\mathbb{P}\circ(u^\varepsilon,v^\varepsilon)^{-1}$, 
we have that for $\tilde{P}-a.s. \omega\in\tilde{\Omega}$ and any $\sigma>0$,
$$
\lim_{\varepsilon\to 0}\left(\int_0^t\int_\mathcal{D}\partial_t\psi(\tau,x)(1+\tilde{u}^\varepsilon(\tau,x))dxd\tau+\int_\mathcal{D}(1+\tilde{u}^\varepsilon(0,x))\psi(0,x)dx-\int_0^t\int_\mathcal{D}\nabla \tilde{v}^\varepsilon\nabla\psi dxd\tau\right)=0,
$$
which yields that
$$
\int_0^t\int_\mathcal{D}\partial_t\psi(\tau,x)(1+\tilde{u}(\tau,x))dxd\tau+\int_\mathcal{D}(1+\tilde{u}(0,x))\psi(0,x)dx-\int_0^t\int_\mathcal{D}\nabla \tilde{v}\nabla\psi dxd\tau=0.
$$
Thus we obtain (\ref{2.8}). The rest proof is the same as the proof of Theorem \ref{t2.3} in Subsection \ref{ss3.5}.

$\hfill\Box$
\vskip.10in

Moreover in radial symmetric case,
\begin{theorem}\label{t5.4.0}
Let $\sigma\geq 0$,
with the same notations as in Theorem \ref{t6.3},
and suppose that the assumptions in Theorem \ref{t6.3} hold. Then in radially symmetric case, that is
$\mathcal{D}=B_1$, where $B_1$ is the unit ball in $\mathbb{R}^d$ and that $u_0^\varepsilon$ is radially symmetric,  we have that
$$
d\mu=2S|D\chi_{E_t}|dxdt\;\;\text{as Radon measure on}\;\;\mathcal{D}_T.
$$
In particular, for $a.e.t\in[0,T]$, $V^t$ is a $(d-1)$-rectifiable varifold (see \cite[Section 11, Section 38]{Simon:1983ti} for definition), i.e.
$$
dV(t,x,p)=2S|D\chi_{E_t}|dxdt\delta_{\vec{\nu}_{E_t}(t,x)}(dp)\;\;\text{as Radon measure on}\;\;\mathcal{D}_T\times P.
$$
Then we have that
\begin{equation}
\left\{
\begin{aligned}
\left(d\mu_{ij}\right)_{d\times d}&=\vec{\nu}_{E_t}\otimes\vec{\nu}_{E_t}d\mu\;\;\text{as Radon measure on}\;\;\mathcal{\bar{D}}_T,\\
v(t,x)&=S\vec{\nu}_{E_t}(x)\cdot \vec{H}_{V^t}(x)\;\;\text{on}\;\;\mathrm{supp}(|D\chi_{E_t}|)\;\;\text{for}\;\;a.e.\;t\in[0,T],
\end{aligned}
\right.
\end{equation}
$\vec{H}_{V^t}$ is the mean curvature vector of $V^t$ defined in Definition \ref{d2.5}
and $\delta_{\vec{\nu}}$ is the Dirac measure concentrated at $\vec{\nu}\in P$.
\end{theorem}

\proof
It suffice to prove
$$
d\mu=2S|D\chi_{E_t}|dxdt\;\;\text{as Radon measure on}\;\;\mathcal{D}_T,
$$
then the following is the same as the proof of Theorem \ref{t2.8} in Section \ref{s4}.

In fact, by taking $h^\varepsilon(t,x)=\int_0^t\int_\mathcal{D}v^\varepsilon(s,x)\xi^\varepsilon(s,x)dsdx$ in (\ref{4.7}), then all the proof followed as in the proof of Theorem \ref{t4.2}. Thus we can finish the proof.

$\hfill\Box$
\vskip.10in

\begin{remark}
The same as in Remark \ref{r2.9}, in radial symmetric case, $v=SH$ on $\Gamma_t$ in a weak sense.  Thus in radial symmetric case, for all $\sigma>0$ the sharp interface limit of equation \ref{6.1} satisfies the deterministic Hele-Shaw model (\ref{1.3}) in a weak sense. In general we also conjuecture that the sharp interface limit of equation \ref{6.1} satisfies the deterministic Hele-Shaw model (\ref{1.3})
\begin{equation*}
\left\{
   \begin{aligned}
   \Delta v&=0 \;\text{in}\;\mathcal{D}\setminus\Gamma_t,\; t>0,\\
\frac{\partial v}{\partial n}&=0\;\text{on}\;\partial\mathcal{D},\\
   v&=SH\;\text{on}\;\Gamma_t,\\
  \mathcal{V}&=\frac{1}{2}(\partial_{n}v^+-\partial_{n}v^-)\;\text{on}\;\Gamma_t.\\
   \end{aligned}
   \right.
\end{equation*}
\end{remark}

Now we will focus on the case that $\sigma=0$. Note that the triple $(E,v,V)$ obtained in Theorem \ref{t6.3} satisfies all the definition in Definition \ref{d2.1} except (\ref{2.8}) and (\ref{2.10}).
Let $\mathcal{D}^+=E_t^o\cap \mathcal{D}$ be the interior of $E_t$ in $\mathcal{D}$ and $\mathcal{D}^-=\mathcal{D}\setminus\bar{E_t}$. 
\begin{theorem}\label{t6.4} Let $(\tilde{\Omega},\tilde{\mathcal{F}},\tilde{\mathbb{P}})$, $E$, $v$ be as in Theorem \ref{t6.3} and $Q$ be an operator safisfying (\ref{2.1}),(\ref{2.2}). Then there exists a $Q$-Wiener process $\tilde{W}$ on $L^2(\mathcal{D})$, which is defined on $(\tilde{\Omega},\tilde{\mathcal{F}},\tilde{\mathbb{P}})$, such that
$$
2d\chi_{E_t}=\Delta vdt+d\tilde{W}_t,
$$
in the sense that
for any $t\in[0,T]$ and 
$\psi\in C_c^1([0,t)\times\bar{\mathcal{D}})$,
\begin{equation}\label{6.10}
\int_0^t\int_\mathcal{D}\left(-2\chi_{E_\tau}\partial_t\psi+\nabla v\cdot\nabla \psi\right)dxd\tau=\int_\mathcal{D}2\chi_{E_0}\psi(0,\cdot)+\int_0^t\langle\psi(\tau,\cdot),d\tilde{W}_\tau\rangle.
\end{equation}
\end{theorem}
\proof  For any $h\in H^1$, denote
$$
M_h^\varepsilon:=\int_\mathcal{D}({u}^\varepsilon(t)-{u}^\varepsilon(0))hdx+\int_0^t\nabla{v}^\varepsilon\cdot\nabla hdx
$$
and
$$
\tilde{M}_h^\varepsilon:=\int_\mathcal{D}(\tilde{u}^\varepsilon(t)-\tilde{u}^\varepsilon(0))hdx+\int_0^t\nabla\tilde{v}^\varepsilon\cdot\nabla hdx.
$$
Clearly,
$$
M_h^\varepsilon=\int_\mathcal{D}h(x)W^\varepsilon_t(x)dx
$$
and as $\varepsilon\to 0$, $M_h^\varepsilon$ converges to a Wiener process with covariance $\|Q^{\frac{1}{2}}h\|_{L^2}^2$. Since $\mathbb{P}\circ(M_h^\varepsilon)^{-1}=\tilde{P}\circ(\tilde{M}_h^\varepsilon)^{-1}$, we know that the law of $\tilde{M}_h^\varepsilon$ converges to a Wiener process with covariance $\|Q^{\frac{1}{2}}h\|_{L^2}^2$.  Moreover
$$
\lim_{\varepsilon\to 0}\tilde{M}_h^\varepsilon=\int_\mathcal{D}({u}(t)-{u}(0))hdx+\int_0^t\nabla{v}\cdot\nabla hdx,\;\;\tilde{P}-a.s..
$$
Thus we obtain that
$$
\int_\mathcal{D}({u}(t)-{u}(0))hdx+\int_0^t\nabla{v}\cdot\nabla hdx
$$
is a Wiener process with covariance $\|Q^{\frac{1}{2}}h\|_{L^2}^2$ on $(\tilde{\Omega},\tilde{\mathcal{F}},\tilde{\mathbb{P}})$. Then there exists a $Q$-Wiener process $\tilde{W}$ on $L^2$, which is defned on $(\tilde{\Omega},\tilde{\mathcal{F}},\tilde{\mathbb{P}})$, such that
$$
\langle W_t,h\rangle=\int_\mathcal{D}({u}(t)-{u}(0))hdx+\int_0^t\int_\mathcal{D}\nabla{v}\cdot\nabla hdxds.
$$

Thus we obtain the following equation for $u$:
$$
du=\Delta vdt+dW_t
$$

Similar to the proof in Subsection \ref{ss3.5},  the It\^o's formula yields that for any $\psi\in C_c^1([0,t)\times\bar{\mathcal{D}})$
\begin{align*}
-\int_\mathcal{D}(1+u(0,x))\psi(0,x)dx=&\int_0^t\int_\mathcal{D}\partial_t\psi(\tau,x)(1+u(\tau,x))dxd\tau-\int_0^t\int_\mathcal{D}\nabla v\nabla\psi dxd\tau\\
&+\int_0^t\langle\psi(\tau,\cdot)dW_\tau\rangle,
\end{align*}
i.e.
$$
\int_0^t\int_\mathcal{D}\left(-2\chi_{E_\tau}\partial_t\psi+\nabla v\cdot\nabla \psi\right)dxd\tau=\int_\mathcal{D}2\chi_{E_0}\psi(0,\cdot)+\int_0^t\langle\psi(\tau,\cdot),d\tilde{W}_\tau\rangle.
$$
Similarly to the disscussion in Subsection \ref{ss2.4}, (\ref{6.10}) is a weak formula for 
$$
2d\chi_{E_t}=\Delta vdt+d\tilde{W}_t.
$$
$\hfill\Box$
\vskip.10in

\begin{corollary}\label{c6.6}
For any $\psi\in C_c^1([0,t)\times\bar{\mathcal{D}})$, with $\mathrm{supp}\psi\subset\mathcal{D}\setminus\Gamma_t$,
$$
\int_0^t\int_{\mathcal{D}}\nabla v\cdot\nabla \psi d\mathcal{H}^dds=\int_0^t\langle\psi,d\tilde{W}_s\rangle.
$$
This is in fact a weak formula for 
$$
\Delta v=-\frac{dW_t}{dt}\;\;\text{in}\;\;\mathcal{D}\setminus\Gamma_t.
$$
\end{corollary}
\proof
Since $\psi\in C_c^1([0,t)\times\bar{\mathcal{D}})$ and $\mathrm{supp}\psi\subset\mathcal{D}\setminus\Gamma_t$, we know that $ \chi_E\psi=\psi$ and $ \chi_E\partial_t\psi=\partial_t\psi$. Thus
$$
\int_0^t\int_\mathcal{D}\chi_{E_t}\partial_t\psi dx=\int_0^t\int_\mathcal{D}\partial_t\psi dx
=-\int_\mathcal{D}\psi(0,\cdot)dx=-\int_\mathcal{D}\chi_{E_0}\psi(0,\cdot)dx.
$$
Then by (\ref{6.10}), we can finish the proof.

$\hfill\Box$
\vskip.10in

\begin{remark}\label{p6.5}
Similar to the deterministic case, $\Delta v$ and $\frac{\partial v}{\partial n}$ are ill-defined. The equation of $(v,\Gamma)$ be understood in distribution sense. We suppose that $(v,\Gamma)$ is smooth enough such that $\Delta v$ and $\frac{\partial v}{\partial n}$ are well-defined.

We also suppose that $\bar{E}\subset \mathcal{D}$.
Denote $\Gamma_t:=\partial E_t$ and let $\mathcal{D}^+=E_t^o\cap \mathcal{D}$ be the interior of $E_t$ in $\mathcal{D}$ and $\mathcal{D}^-=\mathcal{D}\setminus\bar{E_t}$. 

For any $\psi\in C_c^1([0,t)\times\bar{\mathcal{D}})$, with $\operatorname{supp}\psi(t,\cdot)\cap \Gamma_t=\emptyset$,
\begin{equation}\label{6.12}
\begin{aligned}
\int_0^t\int_{\partial\mathcal{D}}\frac{\partial v}{\partial n}\psi d\mathcal{H}^{d-1}ds
=&\int_0^t\int_{\partial\mathcal{D}^-}\frac{\partial v}{\partial n}\psi d\mathcal{H}^{d-1}ds-\int_0^t\int_{\Gamma_t}\frac{\partial v}{\partial n}\psi d\mathcal{H}^{d-1}ds\\
=&\int_0^t\int_{\partial\mathcal{D}^-}\frac{\partial v}{\partial n}\psi d\mathcal{H}^{d-1}ds\\
=&\int_0^t\int_\mathcal{D^-}\text{div}(\nabla v\psi)d\mathcal{H}^d ds\\
=&\int_0^t\int_\mathcal{D^-}\nabla v\cdot\nabla\psi d\mathcal{H}^dds+\int_0^t\int_\mathcal{D^-}\Delta v \psi d\mathcal{H}^d ds\\
=&\int_0^t\langle\psi,d\tilde{W}_s\rangle+ \int_0^t\int_\mathcal{D^-}\Delta v \psi d\mathcal{H}^d ds\\
=&\int_0^t2\langle\psi,d\chi_{E_s}\rangle=0,
\end{aligned}
\end{equation}
where we used  Corollary \ref{c6.6} in the fifth equality. The last  equality holds because $\operatorname{supp}\psi(t,\cdot)\cap E_t=\emptyset$.

Formally we have that in distribution sense
$$
\frac{\partial v}{\partial n}=0\;\;\text{in}\;\;[0,T]\times\partial\mathcal{D}.
$$

To calculate the velocity of $\Gamma_t$, formally we denote $\hat{v}=v+\Delta^{-1}\dot{\tilde{W}}$, where $\dot{\tilde{W}}$ is the formal derivative $\frac{d\tilde{W}}{dt}$. Then we
 have 
$$
2\partial_t\chi_{E_t}=\Delta \hat{v},
$$
and $\hat{v}=0$ in $[0,T]\times(\mathcal{D}\setminus\Gamma_t)$. For any $\psi\in C_c^1(\bar{\mathcal{D}}_t)$
\begin{equation}\label{6.13}
\begin{aligned}
\int_0^t\int_\mathcal{D}\partial_t\chi_{E_t}\psi d\mathcal{H}^dds&=-\frac{1}{2}\int_0^t\int_\mathcal{D}\nabla \hat{v}\nabla \psi d\mathcal{H}^dds\\
&=-\frac{1}{2}\int_0^t\int_\mathcal{D^+}\nabla \hat{v}\nabla \psi d\mathcal{H}^dds-\frac{1}{2}\int_0^t\int_\mathcal{D^-}\nabla \hat{v}\nabla \psi d\mathcal{H}^dds\\
&=\frac{1}{2}\int_0^t\int_\mathcal{D^+}\text{div}(\nabla \hat{v}\psi) d\mathcal{H}^dds+\frac{1}{2}\int_0^t\int_\mathcal{D^-}\text{div}(\nabla \hat{v}\psi) d\mathcal{H}^dds\\
&=\frac{1}{2}\int_0^t\int_{\Gamma_t}(\partial_n \hat{v}^+-\partial_n \hat{v}^-)\psi d\mathcal{H}^{d-1}ds,
\end{aligned}
\end{equation}
which yieds that in distribution sense
$$
\mathcal{V}=\frac{1}{2}(\partial_n \hat{v}^+-\partial_n \hat{v}^-).
$$
Thus formally we have that
\begin{equation}\label{6.14}
\mathcal{V}dt=\frac{1}{2}\left[\frac{\partial}{\partial n}\right]_{\Gamma_t}(vdt+\Delta^{-1} d\tilde{W}_t),
\end{equation}
Here $\left[\frac{\partial}{\partial n}\right]_{\Gamma_t}$ is defined by 
$$
\left[\frac{\partial}{\partial n}\right]_{\Gamma_t}f=\partial_nf^+-\partial_nf^-,
$$
where $f^+,f^-$ is the restriction of $f$ on $\mathcal{D}^+$, $\mathcal{D}^-$, respectively.

\end{remark}

\begin{remark}\label{r5.9}
For the value of $v$ on $\Gamma_t$, since Theorem \ref{t5.4.0} holds for all $\sigma\geq 0$.
Combining it with Corollary \ref{c6.6},  Remark \ref{p6.5} and (\ref{6.14}), we prove that in radial case, when $\sigma=0$, the sharp interface limit of equaiton (\ref{6.1}) is the formally the stochastic Hele-Shaw model (\ref{1.9}). For general case, we conjucture that the sharp interface limit also satisfies (\ref{1.9}):
\begin{equation*}
\left\{
   \begin{aligned}
   \Delta vdt&=-dW_t \;\text{in}\;\mathcal{D}\setminus\Gamma_t,\; t>0,\\
\frac{\partial v}{\partial n}&=0\;\text{on}\;\partial\mathcal{D},\\
   v&=SH\;\text{on}\;\Gamma_t,\\
 \mathcal{V}dt&=\frac{1}{2}\left[\frac{\partial}{\partial n}\right]_{\Gamma_t}(vdt+\Delta^{-1} d{W}_t).
   \end{aligned}
   \right.
\end{equation*}

\end{remark}

 \appendix
  \renewcommand{\appendixname}{Appendix~\Alph{section}}

  \section{Proof of Theorem \ref{t4.2}}
  
  To prove Theorem \ref{t4.2}, we need a technical lemma:
\begin{lemma}\label{l4.2}(\cite[Lemma 5.4]{Chen:1996bg})
For every small positive constant $\delta>0$, there exists a small positive constant $\varepsilon_0(\delta)$ and a large prositive constant $C(\delta)>0$, such that for every $\varepsilon\in(0,\varepsilon_0(\delta)]$, if $(u^\varepsilon,v^\varepsilon)$ is a pair satisfying the second equation in (\ref{4.1}) and 
$$
\|v^\varepsilon\|_{H^1(B_1)}\leq\delta^{-1},\;\;\int_{B_1}e^\varepsilon(u^\varepsilon)dx\leq\mathcal{E}_0,
$$
then the following hold:

(\rmnum{1}). If $(a,b)\subset(\delta,1]$ is an open interval where $|u^\varepsilon|<1-|\ln\varepsilon|^{-\frac{1}{2}}$, then for $a.e.t\in[0,T]$, $u^\varepsilon$ is strictly monotonic in $(a,b)$ and $|b-a|\leq C(\delta)\varepsilon|\ln\varepsilon|$.

(\rmnum{2}). Denote $A^\varepsilon:=\{r\in[2\delta,1-2\delta]:u^\varepsilon(r)=0\}$, then
$$
\int_{2\delta}^{1-2\delta}r^{d-1}e^\varepsilon(u^\varepsilon)dr-C(\delta)\sqrt{\varepsilon}\leq 2S\sum_{r\in A^\varepsilon}r^{d-1}\leq\int_{2\delta-C(\delta)\varepsilon|\ln \varepsilon|}^{1-2\delta+C(\delta)\varepsilon|\ln \varepsilon|}r^{d-1}e^\varepsilon(u^\varepsilon)dr+C(\delta)\sqrt{\varepsilon}.
$$

(\rmnum{3}). For any $r\in A^\varepsilon$,
$$
\left|v^{\varepsilon}(r)+\operatorname{sgn}\left(u_{r}^{\varepsilon}(r)\right) \frac{S(d-1)}{r}\right| \leq C(\delta) \varepsilon^{1 / 8}.
$$

(\rmnum{4}). If $r_1\neq r_2$ in $A^\varepsilon$, then
$$
|r_1-r_2|\geq \frac{1}{C(\delta)}.
$$
\end{lemma}

Now we begin to prove the Theorem \ref{t4.2}:
  
  We ignore the notation tilde $\;\tilde{}\;$ in Theorem \ref{t3.5} for simplicity.

By (\ref{3.4}), Theorem \ref{t2.8} and Proposition \ref{p3.8}, we know that for all $\sigma\geq\frac{1}{2}$, there exists a $h^\varepsilon\in L^2(\Omega,\mathcal{F},\mathbb{P};L^2(0,T))$ such that
\begin{equation}\label{4.7}
\mu^t_\varepsilon(\mathcal{D})=\mathcal{E}^\varepsilon(t)\leq \mathcal{E}_0+h^\varepsilon(t)\;\;\mathbb{P}-a.s.,
\end{equation}
where $d\mu_\varepsilon^t:=e^\varepsilon(u^\varepsilon)dx$ and for $\mathbb{P}-a.s.\;\omega$, $\{h^\varepsilon(\omega,\cdot)\}_\varepsilon$ is bounded in $L^2(0,T)$.

In the following we fix $\omega$ such that all the assertions in Theorem \ref{t3.5} hold, such that (\ref{4.7}) holds, and such that $\{h^\varepsilon(\omega,\cdot)\}_\varepsilon$ is bounded in $L^2(0,T)$.

The following proof is a modification of the proof of \cite[Theorem 5.3]{Chen:1996bg}. We use a contradiction argument. Since $2S|D\chi_{E_t}|dx\leq d\mu^t$, we assume that there exists $T_0\in(0,T]$, such that
$$
\int_0^{T_0}\int_\mathcal{D} d\mu(t,x)=\int_0^{T_0}\int_\mathcal{D}d\mu^tdt>2S\int_0^{T_0}\int_\mathcal{D}|D\chi_{E_t}|dxdt.
$$

Since $du=d\mu^tdt$ is a Radon measure on $\mathcal{D}_T$, we know that $\lim_{\delta\searrow 0}\int_0^T\mu^t(B_\delta)dt=0$ and $\lim_{\delta\searrow 0}\int_0^T\mu^t(B_1\setminus \bar{B}_{1-\delta})dt=0$. Thus there exists $\delta>0$ such that
\begin{equation}\label{4.8}
\int_{0}^{T_0} \int_{B_{1-2 \delta} \setminus \tilde{B}_{2 \delta}} d \mu \geq 2 S\int_{0}^{T_0} \int_{\mathcal{D}}\left|D \chi_{E_{t}}(x)\right| d x d t+\delta\left(T_0+2 \sqrt{T_0\mathcal{E}_0+C_h(\omega)}+1\right),
\end{equation}
where $C_h(\omega):=\sup_{\varepsilon\in[0,1)}\|h^\varepsilon(\omega,\cdot)\|_{L^2(0,T_0)}^2$. For simplicity we denote $C_{T_0}:= \sqrt{T_0\mathcal{E}_0+C_h(\omega)}$. Since $d\mu^\varepsilon:=d\mu^t_\varepsilon dt:=e^\varepsilon(u^\varepsilon)dxdt\to d\mu$, there exists a large positive integer $J\equiv J(\delta)$ such for all $j\geq J$,
$$
\int_{0}^{T} \int_{B_{1-2 \delta} \setminus \bar{B}_{2 \delta}} d \mu_{t}^{\varepsilon_{j}}(x) d t \geq 2 \sigma \int_{0}^{T} \int_{\mathcal{D}}\left|D \chi_{E_{t}}\right| d x d t+\delta(T_0+2 C_{T_0}).
$$

Denote 
$$
\varphi_\delta(t):=\mu^{\varepsilon_j}_t(B_{1-2 \delta} \setminus \bar{B}_{2 \delta}),\;\;\phi(t):=2S|D\chi_{E_t}|(\mathcal{D}).
$$
we have that
\begin{align*}
\int_0^{T_0}\varphi_\delta(t)-\phi(t)dt
\leq& \int_{\{\varphi_\delta-\phi\geq\delta\}}\varphi_\delta(t)-\phi(t)dt+\int_{\{\varphi_\delta-\phi<\delta\}}\varphi_\delta(t)-\phi(t)dt\\
\leq& \delta T_0+ \int_{\{\varphi_\delta-\phi\geq\delta\}}\varphi_\delta(t)-\phi(t)dt\\
\leq& \delta T_0+\int_{\{\varphi_\delta-\phi\geq\delta\}}\varphi_\delta(t)dt\\
\leq& \delta T_0+\|\varphi_\delta\|_{L^2(0,T_0)}\sqrt{\mathcal{H}^1(\{t\in[0,T_0]:\varphi_\delta(t)-\phi(t)\geq\delta\})}\\
\leq& \delta T_0+C_{T_0}\sqrt{\mathcal{H}^1(\{t\in[0,T_0]:\varphi_\delta(t)-\phi(t)\geq\delta\})}.
\end{align*}
In the last inequality we used (\ref{4.7}) and that 
$$
\varphi_\delta(t)\leq \mu^{\varepsilon_j}_t(\mathcal{D})\leq \mathcal{E}_0+h^{\varepsilon^j}.
$$
By (\ref{4.8}), we obtain that
$$
\mathcal{H}^1(\{t\in[0,T_0]:\varphi_\delta(t)-\phi(t)\geq\delta\})\geq 4\delta^2>0.
$$
Moreover
$$
\mathcal{H}^1(\{t\in[0,T_0]:\|v^\varepsilon\|_{H^1}\leq\delta^{-1}\}=1-\mathcal{H}^1(\{t\in[0,T_0]:\|v^\varepsilon\|_{H^1}>\delta^{-1}\}>1-\delta^2\|v^\varepsilon\|_{L^2(0,T;H^1)}^2>0,
$$
where we used that $v^\varepsilon$ converges in $L_w^2(0,T;H^1)$ thus is uniformly bounded in $L^2(0,T;H^1)$.
Hence, for each $j>J$, there exists $t_j\in[0,T_0]$ such that
\begin{equation}\label{4.9}
\|v^{\varepsilon_j}(t_j,\cdot)\|_{H^1}\leq\delta^{-1},\;\;\mu^{\varepsilon_j}_{t_j}(B_{1-2 \delta} \setminus \bar{B}_{2 \delta})\geq 2S|D\chi_{E_{tj}}|(\mathcal{D})+\delta.
\end{equation}
Now we show that (\ref{4.9}) is wrong for $j$ large enough.

For each $j\geq J$, we define
$$
\begin{aligned} 
A^{j} & :=\left\{r \in[\delta, 1-\delta] ; r \in \operatorname{supp}\left(\left|D \chi_{\tilde{E}_{t_{j}}}\right|\right)\right\}, \\ A^{\varepsilon_{j}} & :=\left\{r \in[2 \delta, 1-2 \delta] ; u^{\varepsilon_{j}}\left(r, t_{j}\right)=0\right\} ,
\end{aligned}
$$
where $\tilde{E}\subset[0,T]\times[0,1)$ such that $E=\{(t,x)\in\mathcal{D}_T:(t,|x|)\in\tilde{E}\}$ and $\tilde{E}_t:=\{r\in[0,1):(t,r)\in[0,T]\times[0,1)\}$ is a BV set in $[0,1)$ for any $t\in[0,T]$. By \cite[Theorem 14.3]{Simon:1983ti}, 
\begin{equation}\label{E}
|D\chi_{\tilde{E}}|=\mathcal{H}^0\lfloor \operatorname{supp}
\left(|D\chi_{\tilde{E}}|\right)\;\;|D\chi_{{E}}|=\mathcal{H}^{d-1}\lfloor \operatorname{supp}
\left(|D\chi_{{E}}|\right).
\end{equation}
Since $\mathcal{H}^{0}$ is just the Dirac measure on points in $[0,1)$, we have that $\operatorname{supp}
\left(|D\chi_{\tilde{E}}|\right)=\{r\in[0,1):|D\chi_{\tilde{E}}|(\{r\})\neq 0\}$ which is at most countable. Thus $A^j$ is a countable set. Denote $R(x):=|x|$ is a Lipshitz function on $\mathcal{D}$, thus $R(E)=\tilde{E}$ and $R\left( \operatorname{supp}
\left(|D\chi_{{E}}|\right)\right)=\operatorname{supp}\left(|D\chi_{\tilde{E}}|\right)$. By the area formula \cite[12.6]{Simon:1983ti}, we have that
$$
|D\chi_{E_t}|(\mathcal{D})=\int_0^1\mathcal{H}^{d-1}(\mathcal{D}\cap R^{-1}(y))dy=\int_0^1\omega_dy^{d-1}\chi_{\operatorname{supp}\left(|D\chi_{\tilde{E}}|\right)}(y)dy=\sum_{y\in\operatorname{supp}\left(|D\chi_{\tilde{E}}|\right)}\omega_dy^{d-1}.
$$
Then we have that 
$$
|D\chi_{E_{t_j}}|(\mathcal{D})\geq\sum_{r\in A^j}\omega_dr^{d-1}\geq \delta^{d-1}\omega_d(\#A^j),
$$
where $\#A^j$ is the number of elements in $A^j$ which is finite since $|D\chi_{E_{t_j}}|(\mathcal{D})$ is finite. By the second estimate in Lemma \ref{l4.2}, $A^{\varepsilon_j}$ is also a finite set.

Moreover, by the first inequality in (\ref{4.9}) and Lemma \ref{l4.2}, we have that
$$
\mu_{t}^{\varepsilon_{j}}\left(B_{1-2 \delta} \setminus \bar{B}_{2 \delta}\right)=\int_{2\delta}^{1-2\delta}r^{d-1}e^\varepsilon(u^\varepsilon)dr\leq 2S\sum_{r\in A^\varepsilon}r^{d-1}+C(\delta)\sqrt{\varepsilon}.
$$
Thus since $\omega_d>1$,
 there exists a large integer $J_1\geq J$ such that
$$
\mu_{t_{j}}^{\varepsilon_{j}}\left(B_{1-2 \delta} \setminus \bar{B}_{2 \delta}\right) \leq 2 S \sum_{r \in A^{e_{j}}} \omega_{d} r^{d-1}+\frac{\delta}{2} \quad \forall j \geq J_{1}.
$$
Hence by the second inequality in (\ref{4.8}),
\begin{equation}\label{4.10}
\sum_{r \in A^{\epsilon_j} } \omega_{d} r^{d-1} \geq \sum_{r \in A^{j}} \omega_{d} r^{d-1}+\frac{\delta}{4 S}, \quad \forall j \geq J_{1}
\end{equation}

Denote
$$
l_j:=\sqrt{\varepsilon_j}+\sup_{t\in[0,T_0]}\int_\delta^1|u^{\varepsilon_j}(t,r)+1-2\chi_{\tilde{E}_t}(r)|r^{d-1}dr.
$$
Since $u^{\varepsilon_j}\to -1+2\chi_{E}$ in $C^\beta([0,T];L^1)$,
$$
\sup_{t\in[0,T_0]}\int_\mathcal{D}|u^{\varepsilon_j}(t,x)+1-2\chi_{{E}_t}(x)|dx=\sup_{t\in[0,T_0]}\int_0^1r^{d-1}|u^{\varepsilon_j}(t,r)+1-2\chi_{\tilde{E}_t}(r)|dr
\to 0,$$
for a fixed $\delta>0$,  we have that $\lim_{j\to\infty}l_j=0$.

We claim that the definition of $l_j$ and (\ref{3.10}) imply the existence of $J_2\geq J_1$ such that
\begin{equation}\label{4.11}
\min _{r_{1}, r_{2} \in A^{i}, r_{1} \neq r_{2}}\left|r_{1}-r_{2}\right| \leq 4 l_{j}, \quad \forall j \geq J_{2},
\end{equation}
which is a contradiction to Lemma \ref{l4.2}. 
We prove the (\ref{4.11}) in the following two steps.

First, if $A^{\varepsilon_{j}} \subset \cup_{r \in A^{j}}\left(r-2 l_{j}, r+2 l_{j}\right)$, we claim that  for some $r\in A^j$, there exist at least two elements of $A^{\varepsilon_j}$ in $\left(r-2 l_{j}, r+2 l_{j}\right)$, which concludes (\ref{4.11}). 

If this claim is not true, that is, for any $r\in A^j$, there exista at most one $r_0\in A^{\varepsilon_j}\cap \left(r-2 l_{j}, r+2 l_{j}\right)$. Denote
$$
\underline{A}^j:=\{r\in A^{j}:\exists r_0\in A^{\varepsilon_j}, r_0\in \left(r-2 l_{j}, r+2 l_{j}\right) \}.
$$
Then $\#\underline{A}^j=\#A^{\varepsilon_j}\leq\#A^j$.
Note that
the number of elements in $A^j$ in bounded in $j$ since
$$
\#A^j\leq \delta^{1-d}\omega_\delta^{-1}|D\chi_{E_{t_j}}|(\mathcal{D})= \delta^{1-d}\omega_\delta^{-1}\mathcal{E}^{\varepsilon_j}(t_j)
$$
and $\mathcal{E}^{\varepsilon}(t)$ is uniformly in $\varepsilon$ bounded in $L^\infty(0,T)$.
By (\ref{4.10}) we have that
\begin{align*}
\frac{\delta}{4 S}
=&\sum_{r \in A^{\varepsilon_j}} \omega_{d} r^{d-1}
-\sum_{r \in A^{j}} \omega_{d} r^{d-1}\\
=&\sum_{r \in A^{\varepsilon_j}} \omega_{d} r^{d-1}-\sum_{r \in \underline{A}^{j}} \omega_{d} r^{d-1}
-\sum_{r \in A^{j}\setminus\underline{A}^{j}} \omega_{d} r^{d-1}\\
\leq&\sum_{r \in \underline{A}^{j}} \omega_{d}\left((r+2l_j)^{d-1}-r^{d-1}\right)-\sum_{r \in A^{j}\setminus\underline{A}^{j}} \omega_{d} r^{d-1}\\
\leq& 2l_j\#\underline{A}^j-\sum_{r \in A^{j}\setminus\underline{A}^{j}} \omega_{d} r^{d-1}\\
\leq& C(\delta)l_j-\sum_{r \in A^{j}\setminus\underline{A}^{j}} \omega_{d} r^{d-1},
\end{align*}
which is imposiible for big $j$ since $\lim_{j\to \infty}l_j=0$.

Then if $A^{\varepsilon_{j}} \subset \cup_{r \in A^{j}}\left(r-2 l_{j}, r+2 l_{j}\right)$ does not hold, then there exists $r_1\in A^{\varepsilon_j}$ such that $r_{1} \notin \cup_{r \in A^{j}}\left(r-2 h_{j}, r+2 h_{j}\right)$, i.e. $\left(r_1-2 l_{j}, r_2+2 l_{j}\right)\cap A^j=\varnothing$. Therefore, $\lim_{\varepsilon\to 0}u^\varepsilon\equiv 1$ or $\lim_{\varepsilon\to 0}u^\varepsilon\equiv -1$ on $\left(r_1-2 l_{j}, r_2+2 l_{j}\right)$. Without loss of generality, we assume $\lim_{\varepsilon\to 0}u^\varepsilon\equiv -1$ on $\left(r_1-2 l_{j}, r_2+2 l_{j}\right)$. Thus there exists $(a_1,b_1)\subset(\delta,1)$ such that $r_1\in (a_1,b_1)$ and $|u^{\varepsilon_j}|< 1-|\ln\varepsilon|^{-\frac{1}{2}}$ on $(a_1,b_1)$. By the first assertion of Lemma \ref{l4.2}, $u^{\varepsilon_j}$ is monotonic on $(a_1,b_1)$ and $|b_1-a_1|\leq C(\delta)\varepsilon_j|\ln\varepsilon_j$. Let $\varepsilon_j$ be small enough such that $(a_1,b_1)\subset\left(r_1-2 l_{j}, r_2+2 l_{j}\right)$.

We assume $\partial_ru^{\varepsilon_j}(r_0)>0$, i.e. $u^{\varepsilon_j}$ is monotone increasing on $(a_1,b_1)$. Since $\lim_{\varepsilon\to 0}u^\varepsilon\equiv -1$ on $\left(r_1-2 l_{j}, r_2+2 l_{j}\right)$, we have that
$$
\mathcal{H}^d(u^\varepsilon>0,|u^\varepsilon|>\delta)\leq\int_{B_1\setminus\bar{B}_\delta}|u^{\varepsilon_j}(t,x)+1-2\chi_{{E}_t}(x)|dx\leq l_j.
$$
Since $u^\varepsilon$ is continuous, there must be a $r_2\in(r_1,r_1+2l_j)\cap A^{\varepsilon_j}$. 

In the case that $\partial_ru^{\varepsilon_j}(r_0)>0$, a similar argument yields that there exists $r_2\in(r_1-2l_j,r_1)\cap A^{\varepsilon_j}$. Anyway, we obtain $r_2\in A^{\varepsilon_j}$, such that $|r_1-r_2|\leq 4l_j$. 

Thus we proved (\ref{4.11}), which is a contradiction to Lemma 4.2. Then we finish the proof of the Theorem.

\bibliographystyle{alpha}%

\end{document}